


\magnification 1200
\hsize 13.2cm
\vsize 20cm
\parskip 3pt plus 1pt
\parindent 5mm

\def\\{\hfil\break}


\font\seventeenbf=cmbx10 at 17.28pt

\font\twelvebf=cmbx10 at 12pt
\font\eightbf=cmbx8
\font\sixbf=cmbx6

\font\eighti=cmmi8
\font\sixi=cmmi6

\font\eightrm=cmr8
\font\sixrm=cmr6

\font\eightsy=cmsy8
\font\sixsy=cmsy6

\font\eightit=cmti8
\font\eighttt=cmtt8
\font\eightsl=cmsl8

\font\seventeenbsy=cmbsy10 at 17.28pt

\font\twelvebsy=cmbsy10 at 12pt
\font\tenbsy=cmbsy10
\font\eightbsy=cmbsy8
\font\sevenbsy=cmbsy7
\font\sixbsy=cmbsy6
\font\fivebsy=cmbsy5

\font\tenmsa=msam10

\font\sevenmsa=msam7
\font\fivemsa=msam5
\newfam\msafam
  \textfont\msafam=\tenmsa
  \scriptfont\msafam=\sevenmsa
  \scriptscriptfont\msafam=\fivemsa

\font\tenmsb=msbm10
\font\eightmsb=msbm8
\font\sevenmsb=msbm7
\font\fivemsb=msbm5
\newfam\msbfam
  \textfont\msbfam=\tenmsb
  \scriptfont\msbfam=\sevenmsb
  \scriptscriptfont\msbfam=\fivemsb
\def\Bbb{\fam\msbfam\tenmsb}

\font\tenCal=eusm10
\font\sevenCal=eusm7
\font\fiveCal=eusm5
\newfam\Calfam
  \textfont\Calfam=\tenCal
  \scriptfont\Calfam=\sevenCal
  \scriptscriptfont\Calfam=\fiveCal
\def\Cal{\fam\Calfam\tenCal}

\font\teneuf=eusm10
\font\teneuf=eufm10
\font\seveneuf=eufm7
\font\fiveeuf=eufm5
\newfam\euffam
  \textfont\euffam=\teneuf
  \scriptfont\euffam=\seveneuf
  \scriptscriptfont\euffam=\fiveeuf

\font\seventeenbfit=cmmib10 at 17.28pt

\font\twelvebfit=cmmib10 at 12pt
\font\tenbfit=cmmib10
\font\eightbfit=cmmib8
\font\sevenbfit=cmmib7
\font\sixbfit=cmmib6
\font\fivebfit=cmmib5
\newfam\bfitfam
  \textfont\bfitfam=\tenbfit
  \scriptfont\bfitfam=\sevenbfit
  \scriptscriptfont\bfitfam=\fivebfit


\catcode`\@=11
\def\eightpoint{%
  \textfont0=\eightrm \scriptfont0=\sixrm \scriptscriptfont0=\fiverm
  \def\rm{\fam\z@\eightrm}%
  \textfont1=\eighti \scriptfont1=\sixi \scriptscriptfont1=\fivei
  \def\oldstyle{\fam\@ne\eighti}%
  \textfont2=\eightsy \scriptfont2=\sixsy \scriptscriptfont2=\fivesy
  \textfont\itfam=\eightit
  \def\it{\fam\itfam\eightit}%
  \textfont\slfam=\eightsl
  \def\sl{\fam\slfam\eightsl}%
  \textfont\bffam=\eightbf \scriptfont\bffam=\sixbf
  \scriptscriptfont\bffam=\fivebf
  \def\bf{\fam\bffam\eightbf}%
  \textfont\ttfam=\eighttt
  \def\tt{\fam\ttfam\eighttt}%
  \textfont\msbfam=\eightmsb
  \def\Bbb{\fam\msbfam\eightmsb}%
  \abovedisplayskip=9pt plus 2pt minus 6pt
  \abovedisplayshortskip=0pt plus 2pt
  \belowdisplayskip=9pt plus 2pt minus 6pt
  \belowdisplayshortskip=5pt plus 2pt minus 3pt
  \smallskipamount=2pt plus 1pt minus 1pt
  \medskipamount=4pt plus 2pt minus 1pt
  \bigskipamount=9pt plus 3pt minus 3pt
  \normalbaselineskip=9pt
  \setbox\strutbox=\hbox{\vrule height7pt depth2pt width0pt}%
  \let\bigf@ntpc=\eightrm \let\smallf@ntpc=\sixrm
  \normalbaselines\rm}
\catcode`\@=12

\def\eightpointbf{%
 \textfont0=\eightbf   \scriptfont0=\sixbf   \scriptscriptfont0=\fivebf
 \textfont1=\eightbfit \scriptfont1=\sixbfit \scriptscriptfont1=\fivebfit
 \textfont2=\eightbsy  \scriptfont2=\sixbsy  \scriptscriptfont2=\fivebsy
 \eightbf
 \baselineskip=10pt}

\def\tenpointbf{%
 \textfont0=\tenbf   \scriptfont0=\sevenbf   \scriptscriptfont0=\fivebf
 \textfont1=\tenbfit \scriptfont1=\sevenbfit \scriptscriptfont1=\fivebfit
 \textfont2=\tenbsy  \scriptfont2=\sevenbsy  \scriptscriptfont2=\fivebsy
 \tenbf}

\def\twelvepointbf{%
 \textfont0=\twelvebf   \scriptfont0=\eightbf   \scriptscriptfont0=\sixbf
 \textfont1=\twelvebfit \scriptfont1=\eightbfit \scriptscriptfont1=\sixbfit
 \textfont2=\twelvebsy  \scriptfont2=\eightbsy  \scriptscriptfont2=\sixbsy
 \twelvebf
 \baselineskip=14.4pt}

\def\seventeenpointbf{%
 \textfont0=\seventeenbf  \scriptfont0=\twelvebf  \scriptscriptfont0=\eightbf
 \textfont1=\seventeenbfit\scriptfont1=\twelvebfit\scriptscriptfont1=\eightbfit
 \textfont2=\seventeenbsy \scriptfont2=\twelvebsy \scriptscriptfont2=\eightbsy
 \seventeenbf
 \baselineskip=20.736pt}


\newdimen\srdim \srdim=\hsize
\newdimen\irdim \irdim=\hsize
\def\NOSECTREF#1{\noindent\hbox to \srdim{\null\dotfill ???(#1)}}
\def\SECTREF#1{\noindent\hbox to \srdim{\csname REF\romannumeral#1\endcsname}}
\def\INDREF#1{\noindent\hbox to \irdim{\csname IND\romannumeral#1\endcsname}}
\newlinechar=`\^^J
\def\openauxfile{
  \immediate\openin1\jobname.aux
  \ifeof1
  \message{^^JCAUTION\string: you MUST run TeX a second time^^J}
  \let\sectref=\NOSECTREF \let\indref=\NOSECTREF
  \else
  \input \jobname.aux
  \message{^^JCAUTION\string: if the file has just been modified you may
    have to run TeX twice^^J}
  \let\sectref=\SECTREF \let\indref=\INDREF
  \fi
  \message{to get correct page numbers displayed in Contents or Index
    Tables^^J}
  \immediate\openout1=\jobname.aux
  \let\END=\end \def\end{\immediate\closeout1\END}}

\newbox\titlebox   \setbox\titlebox\hbox{\hfil}
\newbox\sectionbox \setbox\sectionbox\hbox{\hfil}
\def\folio{\ifnum\pageno=1 \hfil \else \ifodd\pageno
           \hfil {\eightpoint\copy\sectionbox\kern8mm\number\pageno}\else
           {\eightpoint\number\pageno\kern8mm\copy\titlebox}\hfil \fi\fi}
\footline={\hfil}
\headline={\folio}

\def\titlerunning#1{\setbox\titlebox\hbox{\eightpoint #1}}
\def\title#1{\noindent\hfil$\smash{\hbox{\seventeenpointbf #1}}$\hfil
             \titlerunning{#1}\medskip}

\newcount\numbersection \numbersection=-1
\def\sectionrunning#1{\setbox\sectionbox\hbox{\eightpoint #1}
  \immediate\write1{\string\def \string\REF
      \romannumeral\numbersection \string{%
      \noexpand#1 \string\dotfill \space \number\pageno \string}}}
\def\section#1{%
  \par\vskip0.666cm\penalty -100
  \vbox{\baselineskip=14.4pt\noindent{{\twelvepointbf #1}}}
  \vskip2pt
  \penalty 500
  \advance\numbersection by 1
  \sectionrunning{#1}}

\def\subsection#1{%
  \par\vskip0.5cm\penalty -100
  \vbox{\noindent{{\tenpointbf #1}}}
  \vskip1pt
  \penalty 500}

\newcount\numberindex \numberindex=0
\def\index#1#2{%
  \advance\numberindex by 1
  \immediate\write1{\string\def \string\IND #1%
     \romannumeral\numberindex \string{%
     \noexpand#2 \string\dotfill \space \string\S \number\numbersection,
     p.\string\ \space\number\pageno \string}}}

\newdimen\itemindent \itemindent=\parindent

\def\item#1{\par\noindent\hangindent\itemindent%
            \rlap{#1}\kern\itemindent\ignorespaces}
\def\itemitem#1{\par\noindent\hangindent2\itemindent%
            \kern\itemindent\rlap{#1}\kern\itemindent\ignorespaces}
\def\itemitemitem#1{\par\noindent\hangindent3\itemindent%
            \kern2\itemindent\rlap{#1}\kern\itemindent\ignorespaces}

\long\def\claim#1|#2\endclaim{\par\vskip 5pt\noindent
{\tenpointbf #1.}\ {\it #2}\par\vskip 5pt}

\def\proof{\noindent{\it Proof}}

\def\today{\ifcase\month\or
January\or February\or March\or April\or May\or June\or July\or August\or
September\or October\or November\or December\fi \space\number\day,
\number\year}

\catcode`\@=11
\newcount\@tempcnta \newcount\@tempcntb
\def\timeofday{{%
\@tempcnta=\time \divide\@tempcnta by 60 \@tempcntb=\@tempcnta
\multiply\@tempcntb by -60 \advance\@tempcntb by \time
\ifnum\@tempcntb > 9 \number\@tempcnta:\number\@tempcntb
  \else\number\@tempcnta:0\number\@tempcntb\fi}}
\catcode`\@=12

\def\bibitem#1&#2&#3&#4&%
{\hangindent=0.8cm\hangafter=1
\noindent\rlap{\hbox{\eightpointbf #1}}\kern0.8cm{\rm #2}{\it #3}{\rm #4.}}


\def\bC{{\Bbb C}}

\def\bQ{{\Bbb Q}}
\def\bR{{\Bbb R}}


\def\cC{{\Cal C}}

\def\cI{{\Cal I}}

\def\cO{{\Cal O}}
\def\cR{{\Cal R}}
\def\cS{{\Cal S}}
\def\cT{{\Cal T}}


\def\square{{\hfill \hbox{
\vrule height 1.453ex  width 0.093ex  depth 0ex
\vrule height 1.5ex  width 1.3ex  depth -1.407ex\kern-0.1ex
\vrule height 1.453ex  width 0.093ex  depth 0ex\kern-1.35ex
\vrule height 0.093ex  width 1.3ex  depth 0ex}}}
\def\qed{\kern10pt$\square$}
\def\hexnbr#1{\ifnum#1<10 \number#1\else
 \ifnum#1=10 A\else\ifnum#1=11 B\else\ifnum#1=12 C\else
 \ifnum#1=13 D\else\ifnum#1=14 E\else\ifnum#1=15 F\fi\fi\fi\fi\fi\fi\fi}
\def\msatype{\hexnbr\msafam}
\def\msbtype{\hexnbr\msbfam}
\mathchardef\restriction="3\msatype16   
\mathchardef\boxsquare="3\msatype03
\mathchardef\preccurlyeq="3\msatype34
\mathchardef\compact="3\msatype62
\mathchardef\smallsetminus="2\msbtype72   
\mathchardef\subsetneq="3\msbtype28
\mathchardef\supsetneq="3\msbtype29
\mathchardef\leqslant="3\msatype36   
\mathchardef\geqslant="3\msatype3E   
\mathchardef\stimes="2\msatype02
\mathchardef\ltimes="2\msbtype6E
\mathchardef\rtimes="2\msbtype6F

\def\ddbar{\partial\overline\partial}

\let\ol=\overline

\let\wt=\widetilde
\let\wh=\widehat
\let\text=\hbox
\def\buildo#1^#2{\mathop{#1}\limits^{#2}}
\def\buildu#1_#2{\mathop{#1}\limits_{#2}}
\def\ort{\mathop{\hbox{\kern1pt\vrule width4.0pt height0.4pt depth0pt
    \vrule width0.4pt height6.0pt depth0pt\kern3.5pt}}}

\def\vlra{\mathrel{\smash-}\joinrel\mathrel{\smash-}\joinrel%
    \kern-2pt\longrightarrow}
\def\srelbar{\vrule width0.6ex height0.65ex depth-0.55ex}
\def\merto{\mathrel{\srelbar\kern1.3pt\srelbar\kern1.3pt\srelbar
    \kern1.3pt\srelbar\kern-1ex\raise0.28ex\hbox{${\scriptscriptstyle>}$}}}




\long\def\InsertFig#1 #2 #3 #4\EndFig{\par
\hbox{\hskip #1mm$\vbox to#2mm{\vfil\special{"
(/home/demailly/psinputs/grlib.ps) run
#3}}#4$}}
\long\def\LabelTeX#1 #2 #3\ELTX{\rlap{\kern#1mm\raise#2mm\hbox{#3}}}


\itemindent = 7mm

\title{A Bergman kernel proof of the}
\title{Kawamata subadjunction theorem}

\titlerunning{Subadjunction}

\vskip10pt

\centerline {\tenrm Bo BERNDTSSON and Mihai P\u AUN}



\vskip20pt

\noindent{\bf Abstract. \it {The main purpose of the following article is to give a proof of  Y. Kawamata's celebrated subadjunction theorem in the spirit of our previous work on 
Bergman kernels. 
We will use two main ingredients : an $\displaystyle L^{2\over m}$--extension theorem 
of Ohsawa-Takegoshi type (which is also a new result) and a more complete version of our former
results.
}}

%
\section{\S0 Introduction}
Let $X$ be a projective manifold and
let $\Theta$ be a closed positive current of (1,1)--type on $X$. A quantitative measure of the logarithmic singularities of $\Theta$ is its {\sl critical exponent} introduced 
e.g. in [10] as follows:

$$C(X, \Theta):= \sup \{c\geq 0: \exp (-c\varphi_\Omega) \hbox { $\in$ } L^1(\Omega) \}$$
for all coordinate set $\Omega\subset X$. The function $\varphi_\Omega$ above is a local potential of $\Theta$ 
on the coordinate set $\Omega$
(see the paragraph 3 for the expression of the normalization we use).
We remark that the above notion makes sense,
since two local potentials of the same current differ by a smooth function. 

One can see that the notion of critical exponent of a closed positive $(1,1)$--current 
is the analytic counterpart of the {\sl log canonical threshold} in algebraic geometry.
Indeed, if we have $\Theta= [D]$ where $D$ is some effective $\bQ$--divisor on $X$, then $C(X, \Theta)= 1$ if and only if $D$ is log canonical. Pushing the analogy with the algebraic geometry
a little bit further, one can easily imagine what the notion of {\sl center} of a current $\Theta$ 
with $C(X, \Theta)= 1$ should be : we have to consider the components of the multiplier ideal sheaf 
associated to $\Theta$. However, it is not yet known that
the multiplier ideal sheaf of such a current is strictly contained in the structural sheaf of 
$X$ ; therefore, in this article we will consider exclusively psh functions with analytic singularities (but nevertheless,  it is possible to obtain
partial results concerning the psh functions whose singularities admit accurate enough approximations with analytic singularities).
\medskip

\noindent Concerning the {\sl minimal center} of a $\bQ$-effective and log canonical divisor
$D$ we have the deep subadjunction theorem of Y. Kawamata, stating that the restriction of the canonical
bundle of the ambient manifold twisted with $D$ to the center is linearly equivalent to 
the canonical class of the center plus a closed positive current.

\medskip

\noindent We establish now the general framework for our results :

{\itemindent 4mm

\item {a)} The current $\Theta$
has analytic singularities ;

\item {b)} $\displaystyle C_{x_0}(X,\Theta)= 1$ (this is the local version of the quantity above, see the paragraph 1.3); 

\item {c)} The cohomology class of $\Theta$ is rational. We let $G$ be a $\bQ$--line bundle whose first Chern class contains $\Theta$.

}

\medskip
\noindent In this perspective, the main goal of the present article is to give a proof of the following version 
of the subadjunction theorem by using the convexity properties of the fiberwise Bergman kernels.

\claim 0.1 Main Theorem|Let $\Theta$ be a closed positive $(1,1)$--current on a projective manifold $X$ which satisfies the properties $a), b)$ and $c)$ above.
 Let $W\subset X$ be the minimal center of $(X, \Theta)$ at $x_0$ ; we assume that 
 it is non-singular. Then there exists a closed positive current $\Theta_{W}$ on $W$ such that
$$\Theta_{W}\in \{{K_X+ G}_{|W}- K_W\},  $$
and moreover the current $\Theta_{W}$ is the weak limit of the currents $\Theta_{W}^{(\varepsilon)}$ such that
$$C_{W}(x_0, \Theta_{W}^{(\varepsilon)})\geq 1.\leqno(\star)$$
\hfill\qed

\endclaim

\medskip
\noindent We remark that the cohomology class of the currents $\Theta_{W}^{(\varepsilon)}$
varies as $\varepsilon\to 0$. Also, the assumption ``$W$ non-singular" above is just 
for the convenience, as if the center $W$ is not smooth, one can see from the proof that 
our result apply in fact to its desingularization.

\bigskip
Let us outline now
the main ideas and statements involved in the proof of the result 0.1. 
\medskip
\noindent The crucial result in the proof of the main theorem is the following
complement to our work [4] ; see equally [3], [19], [29], [35] for related results. We 
remark that the  
{\sl qualitative part} of the theorem below can be seen as a consequence of 
E. Viehweg's results on weak positivity of direct images (again, see [35] and 
the references therein, as well as the recent results of A. H\"oring [17]).  

\claim 0.2 Theorem|Let $p: X\to Y$ be a surjective projective map and let 
$L\to X$ be a line bundle endowed with a metric $h_L$ such that:

{\itemindent 4mm
\smallskip
\item {1)} The curvature current of the bundle $(L, h_L)$ is semipositive, i.e. 
$\displaystyle \sqrt {-1}\Theta_{h_L}(L)\geq 0$;

\item {2)} There exist a generic point $z\in Y$ and a section 
$\displaystyle u\in H^0\big(X_z, mK_{X_z}+ L)\big)$
such that 
$$\displaystyle \int_{X_z}[u, u]_{(m, {{\varphi_L}\over {m}})}< \infty.$$

}

\noindent Then the line bundle $mK_{X/Y}+ L$ admits a metric with positive curvature current.
Moreover, this metric is equal to the fiberwise $m$--Bergman kernel metric 
on the generic fibers of $p$.
\endclaim
\medskip In the statement above the notation under the integral sign is the $\displaystyle L^{2\over m}$ norm of our section $u$. The difference between this result and the one we obtain in [4] is that 
here we allow the map $p$ to be singular (in our previous work we have assumed that $p$ is a smooth fibration).  

The relevant technical result in the proof of 0.2 is the next $L^{2/m}$ extension theorem 
which will be the starting point of this article. The set-up is the following: let $\Omega\subset \bC^n$ be a ball of radius $r$ and let $h:\Omega\to \bC$ be a holomorphic function, such that 
$\sup_\Omega|h|\leq 1$; moreover, we assume that the gradient $\partial h$ of $h$ 
is nowhere zero on the set $V:=  (h=0)$. We denote by $\varphi$ a plurisubharmonic function, such that its restriction to $V$ is well-defined (i.e., $\varphi_{|V}\neq -\infty$). Then the Ohsawa-Takegoshi extension theoren states that for any $f$ holomorphic on $V$, there exists a function $F$, holomorphic in all of $\Omega$, such that $F=f$ on $V$, and moreover

$$\int_{\Omega} |F|^{2}\exp (-\varphi)d\lambda \leq C_0\int_{V} |f|^{2}\exp (-\varphi){{d\lambda_V}\over {|\partial h|^2}}. $$
Here, $C_0$ is an absolute constant. Our generalization is the following.

\claim 0.3 Theorem|For any holomorphic 
function
$f: V\to \bC$  with the property that 
$$\int_{V} |f|^{2/m}\exp (-\varphi){{d\lambda_V}\over {|\partial h|^2}}\leq 1, $$
there exist a function $F\in {\cal O}(\Omega)$ such that :
\item{(i)} $F_{|V}= f$ i.e. the function $F$ is an extension of $f$ ;
\smallskip
\item{(ii)} The next $L^{2/m}$ bound holds 
$$\int_{\Omega} |F|^{2/m}\exp (-\varphi)d\lambda \leq C_0\int_{V} |f|^{2/m}\exp (-\varphi){{d\lambda_V}\over {|\partial h|^2}}, $$
where $C_0$ is the same constant as in the Ohsawa-Takegoshi theorem.
\endclaim

Once this result is established, the proof of 0.2 runs as follows: we first restrict ourselves to the Zariski open set $Y_0\subset Y$ which does not contain any critical value of $p$, and such that 
$\forall z\in Y_0$, all the sections of the bundle $\displaystyle mK_{X_z}+ L$ extend locally near $z$ (the existence of such a set $Y_0$ is a consequence of standard semi-continuity 
arguments). Over $Y_0$, we can apply the result in [4] and therefore the $m$--Bergman kernel metric has a psh variation. We then use the  $L^{2/m}$ extension result to estimate our metric from above by a uniform constant. Standard results of pluripotential theory then gives that the metric extends to a semipositive metric on all of $X$.

The main theorem is proved by resolving the singularities of the current $\Theta$ and by using the properties of the metric constructed 
in 0.2. Indeed, after a few standard perturbation arguments, we will have
$$\mu^*(K_X+ G_{|W}-K_W)+ E_{|W}\equiv K_{S/W}+ \Delta_{\wh X}+  R_{\wh X|W}$$
where $\mu:\wh X\to X$ is a modification which maps the hypersurface $S$ to the minimal center $W$, $E$ is effective and $\mu$-contractible, $\Delta_{\wh X}$ is effective
and  has critical exponent greater than 1 and finally $R_{\wh X}$ is effective and its image via $\mu$ does not contain $W$ (where we still denote by $G, \Delta...$ the approximations of the corresponding bundles for our initial data). 

Our strategy is to show that the current associated to the 
Bergman kernel metric of the bundle 
$$m_0(K_{S/W}+ \Delta_{\wh X})$$
is at least as singular as $m_0[E]$.
Intuitively, this has to be so because the said current admits an extension to $\wh X$ provided that we add to it $m_0R_{\wh X}$ plus the inverse image of a large enough positive 
$(1,1)$--form from $X$. The (crucial) extension properties of closed positive currents are established in the second part of our article ;
the main tools involved in our arguments are regularization results and the invariance of plurigenera.

In conclusion, the existence of this particular representative of the 
Chern class of the right hand side in the formula above show that we can withdrew 
$E$ and still have a pseudoeffective line bundle. The currents $\Theta_{W}^{(\varepsilon)}$ are obtained by taking the direct image of the difference.
We would like to stress on the fact that the explicit expression 
of the $m_0$-Bergman kernel metric is essential for the evaluation of the critical exponent of the currents $\Theta_{W}^{(\varepsilon)}$.\hfill \qed
\smallskip
\noindent In a sequel to this paper [5], the second author  will obtain  analogous {\sl subadjunction--type results for differences of closed positive currents}. The techniques used are somehow similar to those  here, but the whole construction of the proof is a bit more complicated.
\hfill\qed

\claim Notations and conventions|{\rm 

\noindent $\bullet$ The 
metrics of line bundles in this  article are allowed to be singular
unless explicitly stated otherwise (the singularities which we admit are such that the local
weights of the metrics are in $L^1$, so that the curvature current of such object is well defined).

\noindent $\bullet$ Given a metric $h_L$ on a line bundle $L$, we will systematically 
denote by $\varphi_L$ its local weight ; if the notation of the bundle is more complicated 
e.g. $L= \Delta_{\wh X}$ then we will denote its weight by $\varphi_\Delta$.

}

\endclaim

\claim Acknowledgments|{\rm Part of the questions addressed in this article were
brought to our attention by L. Ein and R. Lazarsfeld ; we are very grateful 
to them for this. We would also like
to thank M. Popa and E. Viehweg for interesting and stimulating
discussions about various subjects. Finally, we gratefully acknowledge
support from the Mittag-Leffler institute.

}
\endclaim

\medskip
 
 \vskip 10pt

\vfill
\eject

\section {\S 1. Proof of the main result}
\vskip 10pt We start by
establishing some technical tools which will be needed later in our arguments. In the paragraph 1.1 we deal with a local 
$L^{2/m}$ extension result (theorem 0.3), whereas in 1.2  we prove some positivity results concerning the twisted relative canonical bundle (theorem 0.2). The manner in which these technics are combined in order to prove 0.1 is explained in the subsection 1.3 and 1.4 together with comments about the more general statements/expectations.

\subsection {\S1.1 The $L^{2/m}$ version of the Ohsawa-Takegoshi theorem} 

\vskip 5pt 
In this section we  establish an extension theorem with $L^{2/m}$
bounds, 
analogous to the Ohsawa-Takegoshi result; as we have already mentioned, we will need it for the proof of the theorem 0.1. We recall that the setting is the following: let $\Omega\subset \bC^n$ be a ball of radius $r$ and let $h:\Omega\to \bC$ be a holomorphic function, such that 
$\sup_\Omega|h|\leq 1$; moreover, we assume that the gradient $\partial h$ of $h$ 
is nowhere zero on the set $V:=  (h=0)$. We denote by $\varphi$ a plurisubharmonic function, such that its restriction to $V$ is well-defined (i.e., $\varphi_{|V}\not \equiv -\infty$).
\medskip 
\proof $\hskip 1mm${\sl of the theorem 0.3}.
Let $f$ be a  holomorphic 
function
$f: V\to \bC$  with the property that 
$$\int_{V} |f|^{2/m}\exp (-\varphi){{d\lambda_V}\over {|\partial h|^2}}= 1. $$ 
We begin with some reductions. In the first place we can assume that the function $\varphi$
is smooth, and that the functions $h$ (respectively $f$) can be extended in a neighbourhood of $\Omega$ (of $V$ inside $V\cap \ol\Omega$, respectively). Once the result 
is established under these additional assumptions, the general case follows by approximations and standard normal families arguments.

We can then clearly find some holomorphic $F_1$ in $\Omega$ that extends $f$ and satifies
$$
\int_\Omega |F_1|^{2/m}\exp(-\varphi)d\lambda\leq A<\infty.
$$
We then apply the Ohsawa-Takegoshi theorem with weight
$$
\varphi_1=\varphi +(1-1/m)\log|F_1|^2
$$
and obtain a new extension $F_2$ of $f$ satisfying
$$
\int {{|F_2|^2}\over{|F_1|^{2-2/m}}}\exp-\varphi\, d\lambda\leq
C_0\int {{|f|^2}\over{|F_1|^{2-2/m}}}\exp-\varphi \,{{d\lambda_V}\over{|\partial h|^2}}
=C_0.
$$
H\"older's inequality gives that
$$
\int_\Omega |F_2|^{2/m}\exp(-\varphi)\, d\lambda =
\int_\Omega {{|F_2|^{2/m}}\over{|F_1|^{(2-2/m)/m}}}|F_1|^{(2-2/m)/m}\exp(-\varphi) \,d\lambda\leq 
$$
$$
\leq(\int_\Omega {{|F_2|^2}\over{|F_1|^{2-2/m}}}\exp(-\varphi) d\lambda)^{1/m}
(\int_\Omega|F_1|^{2/m}\exp(-\varphi)d\lambda)^{m/(m-1)}
$$
which is smaller than
$$
C_0^{1/m} A^{m/(m-1)}= A(C_0/A)^{1/m}=:A_1.
$$
If $A>C_0$, then $A_1<A$. We can then repeat the same argument with $F_1$ replaced by $F_2$, etc, and get a decreasing sequence of constants $A_k$, such that 
$$A_{k+1}= A_k(C_0/A_k)^{1/m}$$
for $k\geq 0$. It is easy to see that $A_k$ tends to $C_0$. Indeed, if $A_k >r C_0$ for some $r>1$, then $A_k$ would tend to zero by the relation above. This completes the proof.
\hfill \qed

\vskip 10pt
For further use, let us reformulate the above result in terms of differential forms.
We recall the following notation: if $u$ is a holomorphic $m$--canonical form on an open set $\Omega\subset \bC^p$ and $\varphi$ is some psh function on $\Omega$, then we recall the next notation
$$[u, u]_{(m, \varphi)}:= \big({{(-1)}^{{mp^2}\over {2}}}u\wedge \ol u\big)^{1/m}\exp(-\varphi);$$
it is a (singular) measure on $\Omega$. Then we have the next reformulation of the previous result (we use the same notations as at the beginning of the paragraph).

\claim $0.3^\prime$ Theorem|Given $u$  a holomorphic $m$-canonical form on $V$ such that 
$$\int_{V} [u, u]_{(m,\varphi)}<\infty, $$
there exist a holomorphic $m$-canonical form $U$ on $\Omega$ such that :
\item{(i)} $U= u\wedge (dh)^{\otimes m}$ on $V$ i.e. the form $U$ is an extension of $u$.
\smallskip
\item{(ii)} The next $L^{2/m}$ bound holds 
$$\int_{\Omega} [U, U]_{(m, \varphi)} \leq C_0\int_{V} [u, u]_{(m, \varphi)}, $$where $C_0$ is the constant appearing in the Ohsawa-Takegoshi theorem.

\endclaim
\bigskip

\subsection{\S1.2 Positivity properties of the $m$-Bergman kernel metric}

\smallskip

For convenience of the reader, let us start by recalling some facts concerning the
variation of the fiberwise 
$m$--Bergman metric (also called NS metric in [4], [29], or just a pseudo-norm in [18]).

The general set-up is the following: let $p: X\to Y$ be a surjective projective map between
two non-singular manifolds and consider  a line bundle $L\to X$, endowed with a 
metric $h_L$, such that $\sqrt {-1}\Theta_h(L)\geq 0$. To start with, we will assume 
that {\sl $h_L$ is a genuine metric (i.e. non-singular) and that the map $p$ is a smooth fibration. }

With this data, for each positive integer $m$ we can construct a 
metric $\displaystyle h^{(m)}_{X/Y}$ on the twisted relative bundle 
$\displaystyle mK_{X/Y}+ L$ as follows: let us take a vector $\xi$ in the fiber  
$-(mK_{X/Y}+ L)_x$; then we define its norm 
$$
\Vert \xi\Vert ^2: =\sup |\wt \sigma(x)\cdot \xi|^2
$$
the "$\sup$`` being taken over all sections $\sigma$ to $\displaystyle mK_{X_y}+L$ such that
$$
\int_{X_y} [\sigma, \sigma]_{(m, {1\over m}\varphi_L)}\leq 1;
$$
in the above notation $\varphi_L$ denotes as usual the metric of $L$ and $p(x)= y$; also, we have used the 
notation $\wt \sigma$ to denote the section of the bundle $\displaystyle mK_{X/Y}+ L_{|X_y}$ corresponding to
$\sigma$ via the standard identification (see e.g. [4], section 1). 
\vskip 5pt
For further use, we  give now the expression of the local weights of the metric 
$h^{(m)}_{X/Y}$; we denote by $(z^j)$ respectively $(t^i)$ some local coordinates centered at $x$, respectively $y$. Then we have
$$\exp \big(\varphi^{(m)}_{X/Y}(x)\big)= \sup _{\Vert u\Vert _{m, y}\leq 1}\vert u^\prime(x)|^2$$
where $u\in H^0\big(X_y, mK_{X_y}+L\big)$ is a global section and the above norm means
$$\Vert u\Vert_{m, y}^{2/m}:= \int_{X_y}[u, u]_{(m, {1\over m}\varphi)}.$$
Finally, we denote by $\wt u:= u\wedge (p^*dt)^{\otimes m}$ and the above $u^\prime$ is obtained as  $\wt u= u^\prime (dz)^{\otimes m}$.
\medskip
 In our previous article [4], among other things we have established the fact that {\sl the metric $h^{(m)}_{X/Y}$ above has semi-positive curvature current, 
or it is identically $+\infty$}.
The latter situation occurs precisely when there are no global $L$--twisted $m$-canonical forms on the
fibers. 
\medskip
\noindent In fact, it turns out that the above construction has a meaning even if
the metric $h_L$ we start with is allowed to be singular (but we still assume that the map $p$ is a non-singular fibration). We remark that in this case some fibers $X_y$ may be contained in the unbounded locus of $\displaystyle h_L$, i.e.
$\displaystyle h_{L|X_y}\equiv \infty$, but for such $y\in Y$ we adopt the convention that the 
metric $h^{(m)}_{X/Y}$ is identically $+\infty$ as well. As for the fibers in the complement of this set, 
the family of sections we consider in order to define the metric consists in twisted pluricanonical forms whose $m^{th}$ root is $L^2$. 
\medskip
\noindent We recall now the the result we have proved in the section 3 of [4] (see also [29] where a related result is claimed).
\claim Theorem {[4]}|Let $p: X_0\to Y_0$ be a proper projective non-singular fibration, and let 
$L\to X_0$ be a line bundle endowed with a metric $h_L$ such that:

{\itemindent 4mm
\smallskip
\item {1)} The curvature current of the bundle $(L, h_L)$ is positive, i.e. 
$\displaystyle \sqrt {-1}\Theta_{h_L}(L)\geq 0$;

\item {2)} For each $y\in Y_0$, all the sections of the bundle 
$mK_{X_y}+ L$
extend near $y$;

\item {3)} There exist $z\in Y_0$ and a section 
$\displaystyle u\in H^0\big(X_z, mK_{X_z}+ L)\big)$
such that 
$$\displaystyle \int_{X_z}[u, u]_{(m, {{1}\over {m}}\varphi_L)}< \infty.$$

}

\noindent Then the above fiberwise $m$--Bergman kernel metric $\displaystyle h^{(m)}_{X_0/Y_0}$ has semi-positive
curvature current.\hfill\qed
\endclaim
\bigskip
\noindent We prove now the theorem 0.2.

\proof. By the usual semi-continuity arguments, there exist a Zariski open set $Y_0\subset Y$ such that the condition (2) above is fulfilled; by restricting further the set $Y_0$ we can assume that if we denote by $X_0$ the inverse image 
of the set $Y_0$ via the map $p$, then $p:X_0\to Y_0$ is a non-singular fibration.

Then we can use the above result and infer that our metric $h^{(m)}_0$ is explicitly given over the fibers $X_y$, as soon as $y\in Y_0$ and the restriction of $h_L$ to $X_y$ is well defined.
We are going to show now that this metric admits an extension to the whole manifold $X$. The method of proof is borrowed from our previous work: we show that the local weights of the metric
$h^{(m)}_0$ are locally bounded near $X\setminus X_0$ and then standard results in pluripotential theory imply the result.

Let us pick a point $x\in X$, such that $y:= p(x)\in Y_0$; assume that the restriction of $h_L$ to $X_y$ is well defined. Let $u$ be some
global section of the bundle $\displaystyle mK_{X_y}+ L$. Locally near $x$ we consider a coordinates ball $\Omega$ and let $\Omega_y$ be its intersection with the fiber $X_y$. 
Thus on $\Omega_y$ our section $u$ is just the $m$'th tensor power of some $(n, 0)$ form, and let 
$\wt u:= u\wedge p^{\star}(dt)^{\otimes m}$. With respect to the local coordinates $(z^{j})$ 
on $\Omega$ we have $\wt u= u^\prime dz^{\otimes m}$ and as explained before, the local weight 
of the metric $h^{(m)}_0$ near $x$ is given by the supremum of $|u^\prime|$ when $u$ is
normalized by the condition
$$\int_{X_y}[u, u]_{(m, {1\over m}\varphi)}\leq 1.$$

But now we just invoke the theorem $0.3^\prime$, in order to obtain some $m$--form of maximal degree 
$\wt U$ on $\Omega$ such that its restriction on $\Omega_y$ is just $\wt u$, and such that the next inequality hold
$$\int_\Omega \big[\wt U, \wt U\big]_{(m, {1\over m}\varphi)}\leq C_0 
\int_{\Omega_y}[u, u]_{(m, {1\over m}\varphi)}.$$

The mean value inequality applied to the relation above show that we have  
$$|u^\prime(x)|^{2/m}\leq C\leqno (1)$$ 
where moreover the bound $C$ above does not depend at all on the geometry of the fiber $X_y$, but on the ambient manifold $X$. In particular, the metric 
$h^{(m)}_0$ admit an extension to the whole manifold $X$ and its curvature current is semipositive, so our theorem 0.3 is proved. \hfill\qed
\smallskip
\noindent {\bf 1.2.1 Remark.} For further use, let us point out the fact that the constant $``C "$ above
only depends on  the sup norm of the $m^{th}$--root of the metric of the bundle $L$ and the geometry of the ambient manifold $X$ (and {\sl not at all} on the geometry of the fiber $X_y$). For example, one can easily convinced oneself that if we apply the above arguments to a sequence
$pK_{X/Y}+ L_p$, then the constant in question will be uniformly bounded 
with respect to $p$, provided that $\displaystyle {1\over
  p}\varphi_{L_p}$ is bounded from above.
\hfill\qed
\bigskip
\subsection {\S 1.3 Preliminaries}

\medskip
\noindent Let $X$ be a non-singular projective manifold and let $T$ be a closed positive $(1,1)$--current on $X$. On a coordinate ball 
$\Omega\subset X$ there exist a psh function $\varphi$
such that 
$$T_{|\Omega}= {{\sqrt {-1}}\over {2\pi}}\ddbar \varphi_{\Omega};$$
the function $\varphi_{\Omega}$ above is called {\sl local potential} of the current $T$ 
on $\Omega$.

The notion of local critical exponent which was discussed in the introduction of this article is the following.
\claim 1.3.1 Definition|Let $x_0\in X$ be a point; the quantity 
$$C_{x_0}(X, T):= \sup \{c\geq 0: \exp(-c
\varphi_\Omega) \in  L_{loc}^1(\Omega, x_0), \hbox{ for some }\Omega\subset X,  x_0\in \Omega\}$$ 
is called the critical exponent of the pair $(X,\Theta)$ at $x_0$ ; it is a positive 
real number or it equals $+ \infty$.
\endclaim

\noindent We recall the following notion as well (see [13] and the references therein).

\claim 1.3.2 Definition|The current $T$ has {\sl analytic singularities} if there exist a modification
$\mu:\wh X\to X$ such that 
$$\mu^\star T= \sum_{k=1}^N\nu_j[Y_j]+ \Lambda$$
where $(Y_j)$ are irreducible divisors on $\wh X$, the numbers $(\nu_j)$ are positive reals and $\Lambda$ is a smooth, positive 
(1,1)--form.
\endclaim

\noindent In the right hand side of the above equality, we denote by 
$\mu^\star T$ the current whose local potentials are precisely $\varphi_T\circ \mu$.

One can see that the local potential of a current $T$ with analytic singularities can be written as 
$$\varphi_\Omega= \sum_p\lambda_p\log (\sum_j|f_j^{(p)}|^2)\leqno (2)$$
modulo smooth functions, for any $\Omega\subset X$; in the above expression, $(\lambda_p)$  are positive real numbers which may depend on $\Omega$, and
$(f_j^{(p)})$ are holomorphic functions on $\Omega$.
\medskip

\noindent We will give now a few examples of such currents.
\smallskip
\noindent $\bullet$ Let $L\to X$ be a line bundle and let $T$ be the curvature current
associated to a finite number of sections $(u_j)\subset H^0(X, mL)$. Then $T$ has analytic singularities, as one can see by considering a log resolution of the ideal associated to 
$(u_j)$. \hfill\qed
\smallskip
\noindent $\bullet$ More generally, let $T$ be a closed positive current whose local potentials have the following form
$$\varphi_\Omega= \sum_{p= 1}^N\lambda_p\log (\sum_{j= 1}^M|f_j^{(p,\Omega)}|^2)\leqno (3)$$
modulo smooth functions, where $(\lambda_p)$ are {\sl positive rational numbers}.
Then we claim that $T$ has analytic singularities. To verify this claim, let us consider 
the following coherent ideal
$$\cI^\infty (mT)_x:= \{g\in \cO_{(X, x)}: |g|^2\exp (-m\varphi_{\Omega})< \cO(1)  \}
\leqno (4)$$
We remark that the coherence of this ideal is precisely due to the fact that 
the local potentials of the current are given as in (3).
Next, let us consider a finite covering of the manifold $X$ with the coordinate 
charts $(\Omega_\alpha)$ and let $(\theta_\alpha)$ be a partition of unit subordinate to this cover.
For each index $\alpha$, let us consider the generators $(g_{\alpha, j})$ of the restriction of the ideal $\cI^\infty (mT)$ to $\Omega_\alpha$. We define the current 
$$T_1:= {{{\sqrt -1}} \over{ 2\pi m}}\ddbar \log\big(\sum_\alpha\theta_\alpha^2\sum_j|g_{\alpha, j}|^2\big)$$
We remark that $T_1\geq -C\omega$, where $\omega$ is a fixed metric on $X$ and $C$ is a positive constant (see e.g. [9]). 
The next observation is that the singularities of $T_1$ and $T$ are the same, provided that 
$m$ is divisible enough. Indeed, this is a consequence of the fact that the $m\lambda$ power of each holomorphic function 
$f_j^\Omega$ belong to the ideal $\cI^\infty (mT)$. Therefore, if we consider the log resolution of the ideal $\cI (mT_1)$, the inverse image of the current $T$ will have the decomposition
required by the definition.\hfill\qed

\bigskip

\noindent For the rest of this paragraph, we will assume that the current $\Theta$
verify the assumptions in the statement 0.1, namely, it has analytic singularities, rational cohomology class 
and it is normalized such that 
$$C_{x_0}(X, \Theta)= 1\leqno (5)$$ 
(we remark that this supposes implicitly that the Lelong number of 
$\Theta$ at $x_0$ is positive).

We now rewrite the equality (5)
by using a modification of $X$. To this end, we first remark that we can choose a
birational map $\mu: \wh X\to X$ such that:

{\itemindent 7mm

\item {$(P_1)$} We have
$$\mu^\star(\Theta)= \sum_{p=1}^{N}a^{p}_\Theta[Y_{p}]+ \Lambda$$
where $(Y_p)$ are divisors in $\wh X$ ; we assume that all the exceptional divisors 
of $\mu$ are among this family, therefore some of the coefficients $a^{p}_\Theta\geq 0$ may be zero and $\Lambda$
is a non-singular, semi-positive (1, 1)--form on $\wh X$;

\item {$(P_2)$} The divisors $(Y_{p})$ have normal crossing intersections.
}

\medskip

\noindent Indeed, the existence of $\mu$ is guaranteed by the ``analytic singularities"
hypothesis on the current $\Theta$, together with Hironaka desingularization theorem.
\smallskip
Let $K_{\wh X/X}= \sum_qa^q_{\wh X/X}E_q$ be the relative canonical bundle of the 
map $\mu$. Then for any positive real $t$ we have 
$$
\mu^\star \big(K_X+ tG\big)\equiv K_{\wh X}
+ \sum_{q=1}^{N}\big(ta^{q}_\Theta- a^q_{\wh X/X}\big)[Y_q]+ t\Lambda
\leqno (1\dagger)$$
\noindent 
We recall the cohomology class of the current $\Theta$ is 
assumed to be the same as the Chern class of the $\bQ$-bundle 
$G$ above (see the hypothesis (3) on the page 2).
The
significance of the above formula is that the 
first Chern classes of the corresponding $\bR$-bundles coincide.

Let $\cS_0$ be the set of indexes 
$p\in \{1,..., N\}$
such that $\mu (Y_{p})$ contains the point $x_0$.
The above relation coupled with the change of variables formula show that 
$$\displaystyle \max_{p\in \cS_0} (a^{p}_\Theta C_{x_0}(X, \Theta)- a^q_{\wh X/X})= 1;$$
thus the above relation together with the normalization hypothesis (5) show that 
$$\max_{p\in \cS_0} (a^{p}_\Theta- a^q_{\wh X/X})= 1.\leqno(6)$$
We denote by $\cC(x_0, \Theta, \mu)$ the set of irreducible divisors $(Y_p)$ 
for which the maximum above is obtained.

Next, let us introduce the following notion, analog of the classical one in algebraic geometry.

\claim 1.3.3 Definition|The analytic set $W\subset X$ is called a local center of the pair 
$(X, \Theta)$
at $x_0$ if there exist a birational map 
$\mu$ with the properties $(P_1)$ and $(P_2)$ above, and a divisor $S\in \cC(x_0, \Theta, \mu)$ such that $\mu(S)= W$.
\endclaim

Let us consider the multiplier ideal $\cI(\Theta) \subset \cO_X$, that is 
$$\cI(\Theta)_x:= \{f\in \cO_{X, x}: \int_{(X, x)}|f|^2\exp(-\varphi_{{\Theta}})d\lambda< \infty\}$$
We recall that $\cI(\Theta)$ is coherent ; also, by using a log resolution, we see that the local relative centers are located among the components of the associated zero set of $\cI(\Theta)$ which contain $x_0$, together with their intersections. In particular this show that we have a finite number of local centers
at $x_0$.
\medskip

\noindent Among the several centers a current can have, 
the following one will be very interesting for us.

\medskip
\noindent Let $D$ be an effective $\bQ$-divisor such that its critical exponent at $x_0$ is equal to 1
(i.e. it is lc at $x_0$). In the paper [19], Y. Kawamata show that given
the centers $W_1$ and $W_2$ of $(X, D)$ at some point $x_0\in X$, then any irreducible component
$W$ of the intersection $W_1\cap W_2$ such that $x\in W$ is again a local center of $(X, D)$ at $x_0$.
In this way he shows that there is a well-defined notion of {\sl minimal center } of $(X, D)$
which is unique. 

Our remark is that the same considerations apply in our setting
i.e. the result of Kawamata and Koll\'ar (see e.g. [11] and the references therein) concerning the uniqueness of the minimal center still holds. Indeed, let us give here few details about the technique we use.

\claim Theorem (''Connectedness Lemma`` [11])|Let $\mu: \wh X\to X$ be a 
birational map, where $X$ and $\wh X$ are smooth. Let $\wh D= \sum d_iV_i$ be a $\bQ$-divisor
on $\wh X$, such that $-(K_{\wh X}+ D)$ is $\mu$--nef and big ; we denote the support of the 
multiplier ideal sheaf of the effective part of $\wh D$ with $S$. We assume that if $d_i< 0$, then the corresponding 
hypersurface $V_i$ is contractible.

\noindent Then $S\cap \mu^{-1}(x)$ is connected, for any $x\in X$.

\endclaim 

\noindent We are going to 
use this lemma in our setting as follows. In the first place, we can assume that the
class of the form $\Lambda$ in the inequality 
$(1\dagger)$ is ample, by a slight perturbation of the bundle $G$. Indeed, we recall that 
given any ample 
$A\to X$, the divisor
$$\mu^\star (A)- \sum_j\delta_j E_j$$
is ample on $\wh X$, for small and well chosen positive $(\delta_j)$.

Next, we see that if we denote by 
$$\wh D:= \sum_{p=1}^{N}(a^{p}_\Theta- a^p_{\wh X/X}\big)[Y_{p}]
$$
then we have $-(K_{\wh X}+ \wh D)$ is $\mu$-nef and big,
since it equals $\Lambda$, modulo the pull-back of a line bundle from $X$.
Under these circumstances, one proves that given $W_1$ and $W_2$ two centers,
then any component of their intersection is still a center (see e.g. [19]), thus the notion of {\sl minimal center} has a meaning in our setting.

Let $W$ be the minimal center of 
$(X, \Theta)$ at $x_0$. Together with this data, we have a birational map $\mu: \wh X\to X$ and eventually several divisors with multiplicity 1 dominating $W$ ; the next statement show that we can slightly better.
\claim Claim|
Then there exist a 
divisor $H$ on $X$ and another birational map $\mu_0$ 
such that for any positive $\varepsilon$ the set $\cC(x_0, \Theta, \mu_0)$ 
corresponding to the currents
$\Theta_\varepsilon:= (1-\varepsilon)\Theta+ \delta_\varepsilon H$ 
has an {\sl unique element}. 
\endclaim
 
The argument is absolutely standard (see e.g. [11]), but nevertheless we will give few explanations, for the sake of completeness.

Let $\varepsilon> 0$ be a positive number; we consider the current 
$(1-\varepsilon)\Theta+ \delta_\varepsilon H$, where $H$ is a very ample divisor such that $W\subset H$ and 
such that $H$ do not contain any other center passing thru $x_0$.
We consider a further blow-up of $\wh X$ and call $\mu_0$ the composed map, such that 
the inverse image of $\Theta + [H]$ together with the exceptional divisors have normal crossings;
the equality $(1\dagger)$ become

$$\leqno(2\dagger)
\eqalign{
\mu_0^\star \big(K_X+ (1-\varepsilon)\Theta+ \delta_\varepsilon H\big)=  & K_{\wh X}+
\sum_{p=1}^{N}\big((1-\varepsilon)a^{p}_\Theta- a^p_{\wh X/X}\big)[Y_{p}]+ (1-\varepsilon)\Lambda+\cr
+ &\delta_\varepsilon\big(\sum_{p}a^{p}_H[Y_{p}]+  \Delta\big)\cr
}
$$
where we have denoted the divisors/coefficients involved in the above formula with the same symbols as above.
The last line of $(2\dagger)$ equals $\delta_\varepsilon\mu_0^\star (H)$; we remark that 
the coefficients of $Y_j$ in the first term in the formula above are strictly smaller than 1 : this is a consequence of the fact that the divisorial part of the inverse image of $\Theta$ via $\mu$ 
has normal crossing. 

The support of the divisor $\Delta$ above is disjoint from the other divisors in the decomposition.
If necessary, we can further perturb our data, such that the coefficients 
$a^{p}_H$ are distinct.

Considering all this, one can easily choose $\delta_\varepsilon$ such that the 
previous claim is verified. We remark that the (unique) divisor with coefficient 
equal to 1 will map to the center $W$ precisely 
because of the minimality assumption. In conclusion, the so-called ``tie-break" 
method in algebraic geometry apply in our context without any further modification.
  \hfill\qed

\medskip

\claim Conclusion|In the ``perturbed" setting there exist a unique 
divisor $S$ in $\wh X$ whose  
corresponding coefficient in the formula $(2\dagger)$ is equal to 1, such that 
it dominates the minimal center $W$. \endclaim

\bigskip 

\noindent In order to keep the notations reasonably simple, we still denote by 
$\Theta$ the current obtained after the previous perturbation operations, but we have to keep in mind that at the end of the proof the objects we will get are the $\Theta_{W}^{(\varepsilon)}$
in the statement 0.1. Also, by a slight perturbation of the coefficients $a^{p}$ and 
$b^{(q)}$ we can assume that they are rational (a-priori, they can be arbitrary real numbers, by the definition of {\sl analytic singularities} we adopt here) ; the error term will be absorbed by $\Lambda$. Since the cohomology class of $\Theta$ is rational,
we see that after the perturbation $\{\Lambda\}$ will still be rational.  
\medskip 

Now, we write the relation $(2\dagger)$ as follows:
$$
\mu^\star \big(K_X+ G\big)+ \sum_{p\in \cC_J}(a^p_{\wh X/X}- a^{p}_{\Theta})Y_p  
\equiv K_{\wh X}+ S+ 
\sum_{p\in J}(a^{p}_\Theta-a^p_{\wh X/X})Y_{p}
+ \Lambda.
$$

\noindent In the above relation $J$ denotes the set of integers $p\in \{1,..., N\}$ such that the corresponding coefficient of $Y_p$ is positive ; we denote its complement by $\cC_J$.

\noindent We will denote by $\Delta_{\wh X}$ the $\bQ$--line bundle whose first Chern class contains the following current
$$\sum_{p\in J_0}\big(a^{p}_\Theta- a^p_{\wh X/X}\big)[Y_p]+
 \Lambda\leqno (7)$$
 and such that the relation $(3\dagger)$ below become an equality.
In the definition (7) above we denote by 
  $$J_0:= \{p\in J : x_0\in \mu(Y_q)\}.$$
 
 An important fact is that the critical exponent of the current (7) is greater than 1.
 We equally denote by $R_{\wh X}$ the natural $\bQ$--line bundle associated to the divisor
  $$\sum_{q\in J\setminus J_0}\big(a^{p}_\Theta- a^p_{\wh X/X}\big)
  [Y_p]\leqno (8)$$
and remark that the image via $\mu$ of the support of the current above 
does not contain $W$.

With these notations, the formula above become
$$\mu^*(K_X+ G)+ E= K_{\wh X}+ S+ \Delta_{\wh X}+  R_{\wh X}
\leqno (3\dagger)$$
where
$$E= \sum_{p\in \cC_J}(a^p_{\wh X/X}- a^{p}_\Theta)Y_p $$
is a contractible $\bQ$-divisor on $\wh X$.
\medskip
We recall the next form of the Hartogs extension principle, concerning the line bundles 
$$\wh L:= \mu^\star L+ \wh E$$
where $L\to X$ is a line bundle and $\wh E= \sum b^{(q)}E_q$ is effective and contractible.
The map $\mu: \wh X\to X$ is assumed to be a composition of blow-up maps with non-singular centers. 

\claim Hartog's Lemma|Let $\wh T\in c_1(\wh L)$ be a closed positive current. Then 
$$\wh T\geq [\wh E] ;$$
in particular, the zero set of any holomorphic section of some line bundle numerically equivalent to 
$\wh L$ contains the divisor $\wh E$.
\endclaim

\proof. Let us consider the direct image $T$ of the current $\wh T$ on $X$.
It is a closed positive current in the class $c_1(L)$, and therefore the difference 
$$\wh T- \mu^\star T$$
lies in the class $\{\wh E\}$. On the other hand, we can write
$$\wh T= \chi_{\wh X\setminus E^\prime}\wh T+ \sum \rho^{(j)}E_j$$
and 
$$\mu^\star T= \chi_{\wh X\setminus E^\prime} \mu^\star T+ \sum \tau^{(j)}E_j$$
(where $E^\prime$ is the support of the exceptional divisor).
Since the two currents coincide outside $E^\prime$, we infer that
$$\wh T- \mu^\star T= \sum (\rho^{(j)}- \tau^{(j)})E_j.$$
We invoke now a cohomological argument : the classes associated to the hypersurfaces $E_j$ are independent in $H^{1,1}(X, \bR)$ and therefore
by the above relations we obtain $\rho^{(j)}- \tau^{(j)}= b^{(j)}$ for any $j$,
so our claim is proved.\hfill\qed
  
 \bigskip
The following statement plays an important r\^ole in the proof of the main theorem (see e.g. [11] and the references therein for similar results).

\claim 1.3.5 Proposition|Let $w\in W$ be a generic point and let 
$m\geq 1$ be a positive integer which is divisible enough.
Then the divisor $mE$ is contained in the zero set of any section
$u$ of the 
bundle 
$$m(K_{S_w}+ \Delta_{{\wh X}|S_w}).$$
\endclaim

\smallskip
\noindent {\bf 1.3.6 Remark.} The idea of the following proof is quite simple : 
we first show that the sections above twisted with inverse image of the section of some large enough ample line bundle on $X$ extend to $\wh X$. To this end, the main tools are the theorem 0.2 combined with a version of the Hacon-McKernan {\sl lifting lemma} (see [16])
obtained by Ein-Popa in [14]. Then we apply the Hartogs lemma above.
\hfill \qed

\medskip
\noindent \proof \hskip 1.5mm (of 1.3.5). Since $m$ is assumed to be divisible enough, the multiple $\displaystyle m\Delta_{\wh X}$
is a genuine line bundle. We denote by $p$ the restriction of $\mu$ to $S$.
We want  to apply the theorem 0.2 in order to construct a metric $h^{(m)}_{S/W}$ on the bundle 
$m(K_{S/W}+  \Delta_{\wh X})$ with positive curvature current and such that 
if $w\in W$ is generic, then the restriction of the metric to the fiber $S_w$ has the same singularities
as  the algebraic metric given by the sections of the bundle 
$$\displaystyle m(K_{S_w}+   \Delta_{{\wh X}|S_w}).$$

The metric we use on $\Delta_{\wh X}$ is the one induced by the 
expression (7) and thus we have to be able to show that for 
the generic fiber of the map $p: S\to W$, the space of sections which satisfy the 
$L^{2\over m}$ condition $(2)$ in the theorem 0.2 contains some non-zero element.

Indeed this is the case, since we recall that the critical exponent of $(S, \Delta_{{\wh X}|S})$ 
is strictly greater than 1 (and $w$ is generic), so the corresponding root of 
any non-zero section will be integrable. We remark that we do have such sections, since 
the bundle $mE$ is still effective when restricted to $S_w$, and since $R_{\wh X}$ is vertical with respect to the map $\mu_{|S}= p$.

\noindent Let $A\to X$ be a positive enough line bundle ; we want to apply the Ohsawa-Takegoshi theorem (see the statement 2.2 and the remark 2.3 in the appendix) to the bundle 
$m(K_{S/W}+  \Delta_{\wh X})+\mu^\star A$ in order to extend $u$
(we have already used similar techniques in our article [4]). 

To this end, let us write it as an adjoint bundle as follows
$$m(K_{S/W}+  \Delta_{\wh X})+\mu^\star A= K_S+ (m-1)\big(K_{S/W}+ \Delta_{\wh X}\big)+ \Delta_{\wh X}+ \mu^\star (A-K_W)$$
The second bundle in the above decomposition is endowed with the appropriate multiple of the metric $h^{(m)}_{S/W}$. If we denote by 
$$F:= (m-1)\big(K_{S/W}+ \Delta_{\wh X}\big)+ \Delta_{\wh X}+ \mu^\star (A-K_W)$$
then we have
$$u \in H^0\big(S_w, {K_S+ F}_{|S_w}\big)$$
and let $h_F$ be the metric on $F$ induced by the $h^{(m)}_{S/W}$ together with 
the metric on $\Delta_{\wh X}$ given by the divisors in the decomposition (7) and with the smooth, semi-positively curved metric on 
$\Lambda_2$ and on $\mu^\star(A-K_W)$.
Then we claim that the curvature assumptions and the $L^2$ conditions required by the 
extension theorem are satisfied. 

To verify this claim, let $B\to X$ be a very ample line bundle, such that there exist a family of sections $\rho_j \in H^0(X, B)$ with the property that
$w= \cap_j(\rho_j= 0)$. The curvature conditions to be fulfilled are :

{\itemindent 4mm

\item {i)} $\displaystyle {{\sqrt {-1}}\over {2\pi}}\Theta_{h_F}(F)+{{\sqrt {-1}}\over {2\pi}}\ddbar \log(\sum |\rho_j\circ \mu|^2)\geq 0$
where we measure the norm of the sections $\rho_j$ above by a positively curved metric of $B$;
\smallskip
\item {ii)} $\displaystyle {{\sqrt {-1}}\over {2\pi}}\Theta_{h_F}(F)+ {{\sqrt {-1}}\over {2\pi}}\ddbar \log(\sum |\rho_j\circ \mu|^2)\geq \delta \mu^\star\big({{\sqrt {-1}}\over {2\pi}}\Theta (B)\big)$ where $\delta$ is a positive real.

}
\smallskip
\noindent We see that both conditions will be satisfied if the curvature of $F$
is greater than say $2\mu^\star\Theta (B)$; this last condition can be easily satisfied if $A$ is positive enough.

Concerning the integrability of the section, we remark that we have
$$\int_{S_w}[u, u ]_{({{m-1}\over {m}}\varphi^{(m)}_{S/W}+ \varphi_{\Delta}+\varphi_{A-K_W})}\leq C\int_{S_w}[u, u]_{(m, \varphi_{\Delta}+\varphi_{A-K_W})}<\infty$$
(we refer to the paragraph 1.1 for the meaning of the notation under the integral sign)~;
the first inequality is satisfied because $u$ is bounded with respect to 
$h^{(m)}_{S/W}$ and the second one holds because the ideal sheaf of $\displaystyle \Delta_{\wh X}$ restricted to $S_w$
is trivial, provided that $w$ is generic (which is the case).

Then the Ohsawa-Takegoshi theorem show that $u$ admit an extension 
$U_1 $ as section of the bundle $m(K_{S/W}+ \Delta_{\wh X})+\mu^\star A$. 
We twist the section $U_1$ with the canonical section of the bundle $mR_{\wh X |S}$ and denote by $U$ the resulting section; we remark that it belongs to the
linear system $|m(K_{S/W}+ R_{\wh X}+ \Delta_{\wh X})+ \mu^\star A|$ and that its restriction to the fiber 
$S_w$ is a (non-zero) multiple of our initial section $u$.

In order to simplify the writing, we introduce the next notation
$$R_m:= \Delta_{\wh X}+ R_{\wh X}+ {2\over m}\mu^\star A.$$

\noindent Now let $V\in H^0(W, mK_W+ A)$ which is non-zero at $w$
(again we use the fact that $A$ is positive enough here). Then the section 
$$U\otimes \mu^\star V\in H^0\big(S, m(K_{S}+ R_{m})\big)
\leqno (9)$$
extends (a multiple of) $u$. 

By the definition above, we remark that the bundle $R_{m}$ is naturally decomposed in
two parts:  the $m^{\rm th}$ root of the first one has critical exponent greater than 1
(since it is equal to $\Delta_{\wh X}$) and the section $U\otimes \mu^\star V$ contains 
the second part in the set of its zeroes. In this context, the above section (9) can be extended further to 
$\wh X$ ; the argument is just the Ein-Popa extension theorem [14], see also 
2.4 in the second part of the present article, where we offer (for the convenience of the reader) a proof of their result
``with analytic flavour".
\smallskip

\noindent To see how this result will be applied, the $\bQ$-divisor $\Delta $ in the theorem 2.4 will be simply 
$\displaystyle R_{m}$. The divisors $\Delta_1$ and $\Delta_0$ are respectively the singular part of $R_{\wh X}$ and
$\displaystyle \Delta_{\wh X} $. Finally, the ample part of 
$\displaystyle \Delta_{\wh X}$ plus 
$\displaystyle {2\over m}\mu^\star A$ will be $\Delta_2$.
Next, we recall that we have the relation
$$\mu^*(K_X+ \Theta+ {2\over m}A)+ E= K_{\wh X}+ S+ R_{\wh X}+
\Delta_{\wh X}+ {2\over m}\mu^\star A\leqno (10)$$
The  hypersurfaces $(Y_{p})$ have normal crossing intersections,
and since the supports of $E$ and $S+ R_{\wh X}+ \Delta_{\wh X}$ are disjoint, we see that indeed we can apply the result 2.4 (i.e. the conditions (1)-(4) on the page 29 and
the other hypothesis of 2.4 are fulfilled).

Let 
$$U_{\wh X}\in H^0\big(\wh X, m(K_{\wh X}+ S+ R_{m})\big)$$
be any extension of $U\otimes \mu^\star V$; by the equality $(3\dagger)$, the Chern class of its
zero set  is the same as the one of the
bundle $m\mu^*(K_X+ G+ {2/m}A)+ E$.
By the previous Hartog's lemma, 
the zero set of $U_{\wh X}$ contains the divisor $mE$.
When restricted to $S_w$, we infer that our section $u$ is a multiple of 
$u_E^m$, thus we are done (we use the fact that $V$ {\sl does not} vanish at $w$).
\hfill\qed
\medskip

\bigskip

\subsection {\S 1.4 Proof of the Main Theorem}
\vskip 10pt
\noindent We prove in this subsection the Main Theorem ; as mentioned in the 
introduction, the forthcoming paper [5] will contain a much more 
involved statement and our feeling is that the arguments presented here 
can serve as an introduction to more general facts which will follow. 

The main steps in our proof are the next ones. After few preliminaries, we first show the qualitative part of our statement, i.e. the existence of the closed positive currents 
$\Theta_W^{(\varepsilon)}$ ; this will be a consequence of the theorem 0.2 and 2.1 (in the
second part). We analyze afterwards the size of its critical exponent. 

\medskip

\subsection {\S 1.4.A The singularities of the Bergman metric}
\vskip 5pt 
\noindent To start with, we recall the formula $(3\dagger)$ which states
$$\mu^*(K_X+ G)+ E\equiv K_{\wh X}+ S+ \Delta_{\wh X}+  R_{\wh X}\leqno (12)$$
where $E$ is a contractible, effective $\bQ$ divisor. We recall that the divisor 
$\Delta_{\wh X}$ admits a metric with trivial multiplier sheaf and contains a small ample part (as a result of the perturbation); as for $R_{\wh X}$, it is $\bQ$-effective and trivial on the generic fiber of the map $\mu_{|S}: S\to W$.

We restrict the relation $(12)$ to $S$ and withdraw the inverse image of the canonical bundle of $W$; we denote by $p$ the restriction of the map $\mu$ to $S$ and we get
$$p^*(K_X+ G_{|W}- K_W)+ E\equiv K_{S/W}+ \Delta_{\wh X}+ R_{\wh X}\leqno (13)$$	
(remark that we still denote by $E, \Delta_{\wh X}$ the restrictions of the divisors in questions to $S$). 

Let $m_0$ be a positive integer, such that the multiples $m_0G$, $m_0\Delta_{\wh X}$ and $m_0R_{\wh X}$ become
integral.
We want now to apply our theorem 0.2 in order analyze the positivity properties of the bundle
$$m_0K_{S/W}+ m_0\Delta_{\wh X}.\leqno (14)$$
The dimension of the space of holomorphic sections of its restriction to a generic fiber 
$S_w$ (where $w\in W_0$) is equal to 1, as a consequence of the proposition 1.3.5 ; on the other hand, we can easily identify the space in question : it consists of multiples 
of the divisor 
$m_0E$ restricted to $S_w$. We remark that any section of the above line bundle
automatically satisfies the integrability conditions required by 0.2, since the 
restriction of $\Delta_{\wh X}$ to any generic fiber $S_w$ is still greater than 1.

Thus we have a well-defined $m_0$--Bergman kernel metric $h_{S/W}$
on the bundle (14), whose associated curvature current will be denoted by 
 $\displaystyle {{\sqrt {-1}}\over {2\pi}}\Theta_{S/W}$. We want to show now that this current is {\sl at least as singular as $m_0E$}.
 
To this end, let $A_X\to X$ be a positive line bundle, such that $K_W+ A_X$ is still positive.
We define
$$T_W:= {{\sqrt {-1}}\over {2m_0\pi}}\Theta_{S/W}+ [R_{\wh X}]+ 
\mu^*\big(\Theta(K_W+ A_X)\big)\in c_1\big(K_S+ 
\Delta_{\wh X}+ R_{\wh X}+ \mu^\star(A_X)\big);$$
it is a closed positive current, and our next claim is that there exist a {\sl fixed}
ample bundle $\wh A\to \wh X$ such that for any $m\geq 0$, large and divisible enough, the metric defined by the local potentials of the current above 
on the $\bQ$--bundle 
$K_S+ \Delta_{\wh X}+ R_{\wh X}+ \mu^\star(A_X)$ admits a sub-extension metric $h^{(m)}$ of the 
bundle 
$$\displaystyle K_{\wh X}+ S + \Delta_{\wh X}+ R_{\wh X}+ \mu^\star(A_X)+ 
{1\over m}\wh A,\leqno (15)$$ 
defined on $\wh X$, with positive curvature current. 
The existence of such an object is stated in the theorem 
2.1 in the second part of this article, so let us verify the fact that we are in position to
apply this result.

The $\bQ$-bundle $\Delta$ in the theorem 2.1 will be $\displaystyle R_{\wh X}+ \Delta_{\wh X}$. The divisors $\Delta_1$ and $\Delta_0$ are respectively $R_{\wh X}$ and 
$\displaystyle \Delta_{\wh X} $. If $A_X$ and $m_1$ are large enough, the general member 
of the linear space $H^0\big(\wh X, m_1\mu^*(K_X+ G+ A_X)+ m_1E\big)$
will have normal crossing intersections with the support of the bundles $R_{\wh X}$ and 
$\displaystyle \Delta_{\wh X}$ ; we define $T$ to be the $1/m_1$ times the current of integration over the zero set of such a member. In order to show that with this choice we can apply 2.1 we refer to the part two of our article (pages 30-31). 

Therefore, we get the metric $h^{(m)}$ on the bundle (15) on $\wh X$, with positive curvature current, whose restriction to $S$ is well defined,
and such that 
$$\varphi_{T}\leq \varphi^{(m)}_{|S}+ C(m)\leqno (16)$$
pointwise on $S$. On the other hand, the bundle (15) is numerically equivalent to 
$$\mu^*(K_X+ G+ A_X)+ E+ {1\over m}\wh A$$
and since we can take $\wh A:= \mu^\star B- E_0$ where $B$ is ample
on $X$ and $E_0$ is effective and contractible, the Hartog's lemma
imply that the current associated to $h^{(m)}$ is more singular than $(1-\delta_m)[E]$, where 
$\delta_m\to 0$ as $m\to \infty$ (this little error is induced by $E_0$ above).

Since this is true for all $m$ large enough, we infer from (16) that 
$$T_W\geq [E].\leqno (17)$$
One the other hand the current $T_W$ is a sum of the {\sl positive current} 
$\displaystyle {{\sqrt {-1}}\over {2m_0\pi}}\Theta_{S/W}+ [R_{\wh X}]$ plus a smooth term, therefore we have 
$$ {{\sqrt {-1}}\over {2\pi}}\Theta_{S/W}+ m_0[R_{\wh X}]\geq m_0[E].\leqno (18)$$

As a consequence, the above inequality combined with the formula (13) allows us to define the {\sl closed positive current }
$$\Theta_W\in c_1(K_X+ G_{|W}- K_W)$$ 
on $W$ such that 
$$p^*(\Theta_W)= {{\sqrt {-1}}\over {2m_0\pi}}\Theta_{S/W}+ [R_{\wh X}]-[E] ;$$
this already show that the bundle
$$K_X+ G_{|W}- K_W$$
is pseudoeffective.\hfill \qed
\bigskip

\subsection {\S 1.4.B The critical exponent of $\Theta_W$}
\vskip 5pt 

\noindent  We will show in this paragraph that the critical exponent of the current $\Theta_W$ 
is greater than 1. 

To this end, 
we will analyze during the next few lines the structure of the curvature current of the 
bundle (14), endowed with the metric given in the theorem 0.2.

Let us consider a coordinate set $V\subset W$ ; 
there exist a {\sl meromorphic} section $u_{V}$ of the bundle
$${m_0K_{S/W}+ m_0\Delta_{\wh X}}_{|p^{-1}(V)}\leqno (19)$$
defined by the relation 
$$\displaystyle u_{V}:= {{u_{m_0E}}\over {u_{m_0R_{\wh X}}}}.\leqno(20)$$
One can see that here we implicitly use the fact that the bundle $K_X+ G- K_{W\vert V}$
is trivial, in order to identify the quotient in (20) with a section in the bundle (19). 

We remark that the 
restriction of the meromorphic section above to the {\sl generic fiber} $S_w$ is holomorphic, and it is the only section of the bundle $m_0K_{S/W}+ m_0\Delta_{|S_w}$
up to a multiple. The term ``generic" in the sentence above simply means 
$$w\in W_1:= W\setminus p (R_{\wh X}\cup E^v)$$
In the above relation, we denote by $E^v$ the vertical part (with respect to $p$) of the divisor $E$.
By the discussion at the beginning of the paragraph 1.2, the local weight of the 
metric on $m_0K_{S/W}+ m_0\Delta_{\wh X}$ is given in terms of the section 
$u_{V}$ normalized in a
correct manner, as follows : for $w\in W_1$, let $\tau_w$ be the positive real number such that 
$$\tau_w^2\int_{S_w}{{\vert u_{V}\vert ^{2/m_0}}\over {|\Lambda ^rdp|^2}}
\exp(-\varphi_{\Delta})d\lambda= 1.\leqno (21)$$
One can see that a-priori we need to fix metrics on $S$ and $W$ in order to 
measure the norms above, but they are completely irrelevant, provided that we identify 
$\displaystyle u_{V|S_w}$ with a section of $m_0K_{S_w}+ m_0\Delta_{\wh X|S_w}$, therefore we skip this point, hoping that the confusion caused by it is not too big.

We denote by $\varphi^{(m_0)}_{S/W}$ the local weight of the $m_0$-Bergman kernel metric on the bundle
$m_0(K_{S/W}+ \Delta_{\wh X})$ and  then we have 
$$\exp\big(\varphi^{(m_0)}_{S/W}(x)\big)= \tau_{y}^{m_0}|f_{V}(x)|^2\leqno(22)$$
where $y= p(x)$ and $f_{V}$ is the local expression of the section $u_{V}$ 
near the point $x$ in the fiber $S_y$. It is understood that the equality (22) holds 
for $y\in V\setminus W_1$, and that the content of the theorem 0.2 is that
the expression on the right hand side of (22) is uniformly bounded from above.

\medskip
From the definition of the current $\Theta_W$, the exponential of the $m_0$ times
its
local 
potential is given by the function 
$$x\to \exp\big(\varphi^{(m_0)}_{S/W}\big){{|f_{m_0R_{\wh X}}|^2}\over {|f_{m_0E}|^2}}(x)
\leqno (23)$$
(which in particular, it is log-psh, since $\Theta_W$ is positive). 
On the other hand, the above function is equal to $\tau_{p(x)}$ on $V\setminus W_1$ ;
it is this expression which we will use in order to evaluate the critical exponent of 
$\Theta_W$.

\smallskip
Let $\Omega\subset W$ be a coordinate open set which does not intersect the 
 direct image of $R_{\wh X}$ ;
we want to show that the integral
$$\int _{\Omega}{{d\lambda (w)}\over {\tau^2_w}}$$
converge.

\noindent To this end, remember that by the definition of $\tau$ we have 
$$\tau^2(w)\int_{S_w}{{\vert u_{m_0E}\vert ^{2/m_0}}\over {|\Lambda ^rdp|^2}}\exp(-\varphi_{\Delta}-\varphi_R)= 1.$$
Of course, it may happen that over some fiber 
$S_w$ above the metric of $\Delta_{\wh X}$ is identically $\infty$, or that the $m_0$ root of the section $u_{m_0E}$ does not belong to the multiplier ideal of the restriction of the metric, but this kind of accident can only happen for $w$ in an analytic set;
for such values of $w$ we simply assign the value 0 to $\tau(w)$.

We have 
$$\eqalign{
\int_{\Omega}{{d\lambda (w)}\over {\tau^2w)}}= & \int_{\Omega}d\lambda (w)\int_{S_w}{{\vert u_{m_0E}\vert ^{2/m_0}}\over {|\Lambda ^rdp|^2}}\exp(-\varphi_{\Delta}-\varphi_R)\leq \cr
\leq & \int_{\mu^{-1}(\Omega)}\vert u_{m_0E}\vert ^{2/m_0}\exp(-\varphi_{\Delta}-\varphi_R)<\infty \cr
}$$
where the last inequality is due to the fact that $\Delta_{\wh X}$ has trivial multiplier ideal sheaf. Therefore, the theorem is proved, modulo the fact that we have used a perturbation of $\Theta$ at the beginning. In order to recover the current $\Theta_W$, we simply take the limit in the sense of currents. We would like to mention that 
the multiplier ideal of the limit current may be not trivial. The theorem 0.1 is therefore completely proved.
\hfill \qed

\medskip 
\claim 1.4.7 Remark|{\rm The theorem 0.1 can be seen as a consequence of the original Kawamata's result. Indeed, since the current $\Theta$ has analytic singularities, we can apply e.g. the arguments in [11] to derive the result. 
However, our proof seems to us technically 'lighter'. 
On  the other hand, the information provided by the original approach argument of Kawamata (see equally [11]) is more complete, in the sense that the current $\Theta_W$ is decomposed in a numerically effective part and an effective part in a natural manner.

\hfill\qed}
\endclaim

\medskip 
\claim 1.4.8 Remark|{\rm 
As a consequence of the psh property of the function defined by the relation $(23)$, we see that the normalization 
function $y\to \tau_y$ is {\sl uniformly 
bounded from above}. We recall that this function is only (explicitly) defined on a Zariski open set
and therefore the existence of an upper bound for $\tau$ is a-priori far from obvious,
as one can see by regarding the integral in (21) for $w^\prime$ near some $w\in p(E^v)$.
 Nevertheless it is true, and it imply the next fact : the 
jacobian $\wedge^r dp$ vanish identically on the fiber $S_w\subset E^v$ to a high order, enough to make the 
integral 
$$\int_{S_{w^\prime}}{{\vert u_{V}\vert ^{2/m_0}}\over {|\Lambda ^rdp|^2}}
\exp(-\varphi_{\Delta})d\lambda$$
{\sl bounded from below}, as $w^\prime\to w$. This may be seen as a consequence of the presence of $E$ in the left hand side of (12) : the singularities of some part of the exceptional divisor are bigger than the singularities of the current $\Theta$.
Also, $S$ has a very special status among the exceptional divisors of $\mu$.
\hfill\qed
}
\endclaim

\claim 1.4.9 Remark|{\rm The method of proof we have used 
for the theorem 0.1 can be adapted to a more general context, as follows.

We are giving ourselves  a surjective map
$p: X\to Y$ and a bundle $L\to X$ endowed with a metric with semi-positive curvature, but in addition, let us assume that 
$$\displaystyle h^0\big(X_y, m(K_{X_y}+ L)\big)= 1\leqno (24)$$
for all $m$ large enough and $y\in Y$ generic. Remark that this kind of hypothesis may look odd at a first sight, but 
in fact it turn out to be quite natural in ``subadjunction" context, as the proposition 1.3.5 shows it. Moreover we will assume that the metric 
$h_L$ has trivial multiplier sheaf and analytic singularities. Then we have the next consequence of our previous results.
\medskip 
\claim 1.4.10 Corollary|Under the above hypothesis, there exist metric on the twisted relative canonical bundle 
$$K_{X/Y}+ L$$
whose curvature current equals
$${{\sqrt {-1}}\over {2\pi}}\Theta_{X/Y}:= [\Xi]+ \chi_Zp^*(R_Y)$$
where $\Xi$ is an effective $\bR$-divisor on $X$ such that its intersection to the generic fiber is the current of integration over the zero set of the section in (24) renormalized, $Z$ is $p$-vertical and finally $R_Y$ is a closed positive current on $Y$.
\endclaim
\medskip 
\proof. The argument is derived from the theory of closed positive currents, 
(especially the structure theorems of Y.-T. Siu, see [30]), combined with our theorem 0.2. Indeed, let us 
consider the curvature current $\Theta_{X/Y}$ of the metric $h^{(m)}_{X/Y}$; it is closed, positive, and its restriction to the fiber $X_y$ is given by the current of integration over the analytic set corresponding to the section in $m(K_{X_y}+ L)$ normalized by the factor $1/m$, if $y\in Y\setminus Z$, where $Z$ is a Zariski closed set. 

On the other hand, if $m\gg 0$ then there exist a divisor $D_m$ on $X$ whose intersection with the generic fiber of $p$ is the zero set of the section in the linear system 
$m(K_{X_y}+ L)$; this fact is also a consequence of usual extension theorems, together with the hypothesis (24). Now by the results in [29], the generic Lelong number of $\Theta_{X/Y}$ along $D_m$ is the same as the generic Lelong number 
of the restriction of $\Theta_{X/Y}$ to a generic fiber along the corresponding restriction of $D_m$. Thus we have
$$\chi_{X\setminus Z}{{\sqrt {-1}}\over {2\pi}}\Theta_{X/Y}\geq {1\over m}\chi_{X\setminus Z}[D_m]$$
and the restriction of the difference of the currents above to the generic fibers of $p$ is equal to zero.
Positivity considerations show the existence of a closed positive current $R_Y$ on $Y$ such that 
$${{\sqrt {-1}}\over {2\pi}}\Theta_{X/Y}- {1\over m}D_m= p^\star (R_Y)$$
holds on the complement of the Zariski closed subset of $Z\subset X$: indeed, one first show that such a current exists on 
the complement of the critical values of $p$, and then that is extends since it has bounded mass near this set.
Remark that in this way we have determined the part of the current 
$\displaystyle{{\sqrt {-1}}\over {2\pi}}\Theta_{X/Y}$ outside $Z$.
But the part of the current supported in $Z$ is of divisorial type as well, by the support theorem
[30], thus this gives the current $\Xi$ and the corollary is proved.
\hfill\qed}
\endclaim

\medskip

\claim 1.4.11 Remark|{\rm Let $\Theta$ be a closed positive current on $X$ of
(1,1)--type. We assume that its critical exponent at some point $x_0\in X$ is equal to 1.
Concerning the setting we have the next so-called {\sl openess conjecture}, due to Demailly-Koll\'ar in [10].

\claim Conjecture ([10])|The integral 
$$\int_{(X, x_0)}\exp (-\varphi_T)d\lambda$$
is divergent.
\endclaim

As far as we are aware, the previous statement was only established 
for surfaces (see [15]). It would be extremely useful to have it for arbitrary dimensions, 
since in our opinion this is the missing piece in order to prove the theorem 
0.1 without the {\sl analytic singularities} hypothesis. 
}
\endclaim

\bigskip

\section {\S 2. Extension of metrics}

\medskip

\noindent In this second part we  prove 
an extension theorem concerning  metrics with positive curvature current of adjoint 
$\bQ$--line bundles. Particular cases of the statement below (and its corollaries) were used several times in the 
proof of the main theorem. 

\noindent The general framework of our result is as follows. 

\noindent $\bullet$ We consider a projective, non-singular manifold $X$,
and let $\Delta_0\to X$ be a $\bQ$-line bundle, such that for any $m$ large and divisible enough, we have
$$m\Delta_0\equiv L_1+\cdots L_{m-1}\leqno (25)$$
where for each $1\leq j\leq m-1$ the $L_j$ is a line bundle, endowed with a 
metric $h_j$ whose curvature current is assumed to be positive. 
\smallskip
\noindent $\bullet$ $\Delta_1\to X$ is a $\bQ$--line bundle, whose first Chern class
contains a closed positive current $\Xi$. We introduce the  notation
$$\Delta:= \Delta_0+ \Delta_1.$$
\smallskip

\noindent $\bullet$  $S\subset X$ is a non-singular hypersurface such that $h_{j|S}$ 
and $\Xi_{|S}$ are well-defined,
in the sense that the local potentials of $h_j$ are not identically $-\infty$ when restricted to $S$. We assume the existence of a closed positive current 
$$T\in c_1(K_X+ S+ \Delta)$$
such that :
{\itemindent 6mm

\smallskip
\item {$(\cR)$} The restriction of $T$ to $S$ is well defined ; 
\smallskip
\item {$(\cT)$} There is some positive number $\varepsilon$ such that for any $\delta>0$ and  any $j= 1,..., m-1$ we have 
$$\cI \big( (1-\delta)\varphi_{L_j}+ \varepsilon \varphi_{T|S}\big)= \cO_S.\leqno (\cT_j)$$
}
\medskip

\noindent We observe that under the hypothesis $(\cT_j)$ the multiplier sheaf
of the metric $(1-\delta)\varphi_j$ is automatically trivial for any $\delta> 0$ . The assumption above
intutively means  that the singularities of $\varphi_j$ and $\varphi_T$ are transversal.
\bigskip
The following statement is a {\sl metric version} of the extension theorem due to
Ein-Popa and Hacon-McKernan [14], [16] (as it will be seen in the corollaries which  follow  the proof). The origin of  this kind of results is the seminal article of Y.-T. Siu  [32].

\claim Theorem 2.1|Let $S$ be a non-singular hypersurface of a projective manifold $X$
and let $\Delta_0, \Delta_1, T$ be respectively $\bQ$--line bundles and a closed positive current with the properties stated above.
We consider a metric $h_{K_S+ \Delta}$ on the $\bQ$-line bundle
$K_S+ \Delta$ such that :
{\itemindent 5mm

\smallskip
\item {i)} The corresponding curvature current is positive ;
\smallskip
\item {ii)} We have
$$\varphi_{K_S+ \Delta}\leq \varphi_{\Xi|S},$$ 
or more generally, it is enough to assume that for any positive integer $m$ we have
$$\cI (m\varphi_{K_S+ \Delta})\subset \cI(m\varphi_{\Xi|S}).\leqno (26)$$
}
\noindent Then there exist an ample line bundle $A\to X$ such that 
for any positive 
integer $m$,  
the metric $\displaystyle h_{K_S+ \Delta}$ admits a sub-extension 
$\displaystyle h^{(m)}$, where  $\displaystyle h^{(m)}$ is a metric of the $\bQ$--line bundle
$K_X+ S+ \Delta+1/mA$ defined on $X$, with positive curvature current.\hfill \qed
\endclaim
\medskip
\noindent In other words, we claim the existence of a metric  
$\displaystyle h^{(m)}$ of the $\bQ$--line bundle
$K_X+ S+ \Delta+1/mA$ with positive curvature and such that its restriction to $S$ is well defined and satisfies the next inequality
$$\varphi_{K_S+ \Delta}\leq \varphi^{(m)}- \varphi_A/m,$$ 
at each point of $S$, where $\varphi_A$ is a smooth metric on $A$.

\proof. Our arguments involve two main ingredients : an approximation theorem for  closed positive currents (which is originally due to Demailly in [9], and further polished by Boucksom in [6]) and the invariance of plurigenera  techniques. The first technical tool allows 
us to ``convert" the metric $\displaystyle h_{K_S+ \Delta}$ to a family of sections of the bundle $m(K_S+ \Delta)+ A$, whereas the second will show that each member of the
family of sections obtained in this way admits an extension to $X$. We remark that for the second step we could just quote
the recent result of Ein-Popa [14] ; however, we prefer to present an alternative proof, which is more analytic in spirit (the reader can profitably consult the original proof in [16], based on the 
theory of {\sl adjoint ideals}).
\medskip

\subsection {\S 2.1 Approximation of closed positive currents}
\smallskip
\noindent To start with, we  recall the (crucial) fact that a closed positive current of 
$(1, 1)$--type
can be approximate in a very accurate way by 
a sequence of currents given by algebraic metrics.

\claim Theorem A ([6], [9])|There exist an ample line bundle $A \to X$
and a smooth metric $h_A$ on $A$ 
such that the following holds. Let $m$ be a large and divisible enough positive integer, and let $\big(\sigma_j^{(m)}\big)$ be an orthonormal basis of the space
$$H^0\big(S, m(K_S+ \Delta)+ A\big)$$
endowed with the metric whose local weights are given by the functions
$\displaystyle m\varphi_{K_S+ \Delta}+ \varphi_{A}$ and an arbitrary smooth volume form on $S$. Then we have the pointwise inequality
$$\varphi_A/m+\varphi_{K_S+ \Delta}\leq {1\over m}\log\big(\sum_{j=1}^{N_m}|f_j^{(m)}|^2\big)+ C\leqno (27)
$$
where $f_j^{(m)}$ is the local holomorphic function corresponding to the 
section $\sigma_j^{(m)}$, for any $j$.\hfill \qed

\endclaim

\medskip
\noindent We remark that the precise statement in [6], [9] is slightly different, 
but the proof given there clearly contains the inequality (27) above. 
Notice that if the metric $\varphi_{K_S+ \Delta}$ satisfies condition (ii) of Theorem 2.1, then 
the sections $(\sigma_j^{(m)})$ satisfy the integrability relation 
$$\sigma_j^{(m)}\in \cI (m\varphi_{\Xi |S})\leqno (28)$$
for any positive integers $j, m$.
\medskip

\subsection {\S 2.2 Extension of sections}
\smallskip

\noindent Let $u$ be a section of $m(K_S+\Delta)+A$
such that 
$$\int_S\vert u\vert ^2\exp (-m\varphi_{K_S+ \Delta}-\varphi_A)d\lambda< \infty;\leqno (29)$$
so that in particular $u\in \cI (m\varphi_{\Xi |S})$ (for example, $u$ may be 
one of the $\sigma_j^{(m)}$ involved in the approximation process above).

The main part of the proof of 2.1 is to show that $u$ extends to $X$.
(Of course, we are going to use heavily the other metric assumptions 
in the statement 2.1, i.e. just the relation (29) above is  not sufficient). 
To this end, we will use the invariance of plurigenera techniques. 
The argument is derived from the original article
of Siu, see [32], as well as [13], [26].

\noindent We denote by $L_m:= m\Delta_1$ and
$$L^{(p)}:= p(K_X+ S)+L_1+...+ L_p\leqno (30)$$
where $p= 1,..., m$. By convention, $L^{(0)}$ is the trivial bundle.
We can assume that the bundle $A\to X$ in 
the approximation statement above is ample 
enough, but independent of $m$, so that the next two conditions are satisfied.

{
\itemindent 6mm
\noindent
\item {$(\dagger)$} For each $0 \leq p\leq m-1$, the bundle
$L^{(p)}+ mA$ is generated by its global sections, which we denote by
$(s^{(p)}_j)$;
\smallskip
\item {$(\dagger\dagger)$} Any section of the bundle $L^{(m)}+
mA_{\vert S}$ admits an extension to 
$X$.

}
\medskip
\noindent  
We formulate now the following inductive statement (depending on a positive integer $k$ and 
on $0\leq p\leq m-1$) :
\vskip 7pt
\noindent 
{$\bigl({\cal P}_{k,p}\bigr):$ 
\sl The sections 
$$u^k\otimes s^{(p)}_j\in H^0\bigl(S, L^{(p)}+ kL^{(m)}+ (k+m)A_{\vert S}\bigr)$$
extend to $X$, for each $j= 1,..., N_p.$}
\vskip 7pt
\noindent If we can show that ${\cal P}_{k,p}$ is true for any $k$ and $p$, 
then we are done, since the extensions $U^{(km)}_j$ of the sections 
$u^k\otimes s^{(0)}_j$ can be used to define a metric
$$
\log\sum|U_j^{(km)}|^2
$$
on the bundle
$$km(K_X+S+ \Delta)+ (k+m)A$$
whose $km^{\rm th}$ root it is defined to be 
$h^{(m)}$ (for $k\gg 0$). Theorem A shows that  $h^{(m)}$
is then indeed a sub-extension of the given metric $h_{K_S+ \Delta}$.

\medskip
\noindent Thus, it is enough to check the property ${\cal P}_{k,p}$ ; our
arguments rely heavily on the next version of the
Ohsawa-Takegoshi theorem, proved by L. Manivel (see [1], [2], [12], [22], [24], [25] to quote
only a few), and further refined by McNeal-Varolin in [23].

\claim 2.2 Theorem ([23])|Let $X$ be a projective
$n$-dimensional 
manifold, and let $S\subset X$ be a non-singular hypersurface.
Let $F$ be a line bundle,
endowed with a metric $h_F$. We assume the existence 
of some non-singular metric $h_S$ on the bundle $\cal O(S)$ such that :  

{\itemindent 6mm
\smallskip
\item {(1)} $\displaystyle {{\sqrt {-1}}\over {2\pi}}\Theta_F
\geq 0$ on $X$;
\smallskip
\item {(2)} $\displaystyle {{\sqrt {-1}}\over {2\pi}}\Theta_F -  
\alpha
 {{\sqrt {-1}}\over {2\pi}}\Theta_S\geq 0$ 
for some $\alpha> 0$;
\smallskip
\item {(3)} The restriction
of the metric $h_F$ on $S$ is well defined.

}
\smallskip

\noindent Then every section $u\in H^0\bigl(S, (K_X + S+
F_{\vert S})\otimes {\cal I}(h_{F|S})\bigr)$ admits an extension
$U$ to $X$ such that 
$$c_n\int_XU\wedge \ol U\exp \Big(-\varphi_F-\varphi_S-\log\big(|s|^2\log^2(\vert s\vert)\big)\Big)< \infty$$
where $s$ is a section whose zero set is precisely the hypersurface $S$ and its norm in the integral above is measured with respect to $h_S$. In particular
$$
c_n\int_XU\wedge \ol U\exp \Big(-\varphi_F-\varphi_S\Big)/|s|^{2(1-\delta)}< \infty$$
for any $\delta>0$.
\endclaim
\medskip 
\claim 2.3 Remark|{\rm In fact, the general result in [23] allow to extend sections 
of bundles $K_X+ E_{|W}$ defined on higher codimensional manifolds $W\subset X$
which are complete intersections 
$$W= \cap_j (\rho_j= 0)$$
where $\rho_j$ are sections of some line bundle $\Lambda$. Of course, in this case the curvature and integrability conditions $(1)-(3)$ are modified : basically we replace 
$\displaystyle {{\sqrt {-1}}\over {2\pi}}\Theta_F$ with 
$\displaystyle {{\sqrt {-1}}\over {2\pi}}\Theta_E+ \displaystyle {{\sqrt {-1}}\over {2\pi}}
\ddbar \log \big(\sum |\rho_j|^2\big)$, and the curvature term of $\cO(S)$ with 
the curvature of $\Lambda$  
(if $W= S$ is a hypersurface, the bundle $E$ in this setting 
is just $F+ S$). We refer to [23] for the precise statement.\hfill \qed
}
\endclaim

\medskip
We will use inductively the theorem above, in order to show the extension property
${\cal P}_{k,p}$. The first steps are as follows.

\noindent (1)  For each $j= 1,...,N_0$, the section 
$u\otimes s^{(0)}_j\in H^0\bigl(S, L^{(m)}+ (m+1)A_{\vert S}\bigr)$
admits an extension $U^{(m)}_j\in H^0\bigl(X,  L^{(m)}+ (m+1)A\bigr)$,
by the property $\dagger\dagger$.

\noindent (2) We use the sections $(U^{(m)}_j)$ to construct a metric 
$$\varphi^{(m)}=\log\sum|U_j^{(m)}|^2$$
on the bundle $L^{(m)}+ (m+1)A$.

\noindent (3) Let us consider the section 
$u\otimes s^{(1)}_j\in H^0\bigl(S, L^{(1)}+ L^{(m)}+ (1+m)A_{\vert S}\bigr)$. 
We remark that the bundle
$$L^{(1)}+ L^{(m)}+ (m+1)A = K_X+ S+ L_1+ L^{(m)}+ (m+1)A$$
can be written as $K_X+ S+ F$ where
$$F:=  L_1+ L^{(m)}+ (m+1)A.$$
We are going to construct a metric on $F$ which satisfy the 
curvature and integrability assumptions in the theorem 2.2.

Let $0< \delta\ll \varepsilon\ll 1$ be positive real numbers. We endow the bundle 
$F$ with the metric given by
$$\varphi^{(m)}_{\delta, \varepsilon}:= (1-\delta)\varphi_{L_1}+ \delta\wt \varphi_{L_1}+ 
(1-\varepsilon)\varphi^{(m)}+ \varepsilon((m+1)\varphi_A+ m\varphi_T)\leqno (31)$$
where the metric $\wt \varphi_{L_1}$ is smooth (no curvature requirements) and  
$\varphi_{L_1}$ is the singular metric whose existence is given by hypothesis.
The curvature conditions (1) and (2) in 2.2 are satisfied, since we can use the metric on $A$
in order to absorb the negativity of $\wt \varphi_{L_1}$. 

Next we claim that the sections $u\otimes s^{(1)}_j$ are integrable with respect to the metric defined in (31).
Indeed, the integrand is less singular than 
$$
\exp \big(-(1-\delta)\varphi_{L_1}- \varepsilon\varphi_T\big)$$
which is integrable  by the hypothesis $(\cT_1)$.
\smallskip

\noindent (4) We apply the extension theorem 2.2 and we get
$U^{(m+ 1)}_j$, whose
restriction on $S$ is precisely
$u\otimes s^{(1)}_j$.\hfill \qed
\medskip

\vskip 5pt Now the assertion  ${\cal P}_{k,p}$ will be obtained by iterating
the
procedure (1)-(4) several times. Indeed, assume that the 
proposition ${\cal P}_{k,p}$ has been checked, and consider the set of global sections 
$$U^{(km+p)}_j\in H^0\bigl(X, L^{(p)}+ kL^{(m)}+ (k+m)A\bigr)$$
which extend $u^k\otimes s^{(p)}_j$. They induce a metric on the above bundle, denoted by $\varphi^{(km+p)}$. 

\noindent If $p< m-1$, then we define the family of 
sections 
$$u^k\otimes s^{(p+1)}_j\in H^0(S, L^{(p+1)}+ kL^{(m)}+ (k+m)A_{|S})$$
on $S$. 
As in the step (3) above  we have that
$$L^{(p+1)}= K_X+ S+ L_{p+1}+ L^{(p)}.$$
To apply the extension result 2.2, we have to exhibit a metric
on the bundle 
$$F:=  L_{p+1}+ L^{(p)}+kL^{(m)}+ (k+ m)A$$
for which the curvature conditions are satisfied, and such that the family of sections above are $L^2$ with respect to it.
We define
$$\varphi^{(km+p+1)}_{\delta, \varepsilon}:= (1-\delta)\varphi_{L_{p+1}}+ \delta\wt \varphi_{L_{p+1}}+ 
(1-\varepsilon)\varphi^{(km+ p)}+ \varepsilon( \wt \varphi_{L^{(p)}}+ mk\varphi_T+ (m+k)\varphi_A)\leqno (32)$$
and we see that the parameters $\varepsilon, \delta$ have to satisfy the conditions :

{\itemindent 6mm

\item {(A)} We don't assume any relation between the zero set of $u$ and the singularities of $T$, thus we have to keep the poles of $T$ ``small" in the expression of the metric above, as we have to use 
$(\cT_{p+1}$). This impose $mk\varepsilon\ll 1$ ;

\item {(B)} We have to absorb the negativity in the smooth curvature terms in (32), thus we are forced to have $\delta\ll\varepsilon$.

}
\noindent Clearly we can choose $\varepsilon, \delta$ with the properties required above. Let us check next the $L^2$ condition. As before, the integrand is less singular than
$$
\exp \big(-(1-\delta)\varphi_{L_{p+1}}- mk\varepsilon\varphi_T\big).$$
This is again integrable  because of the transversality hypothesis $(\cT_{p+1}$).

\medskip
Finally, let us indicate how to perform the induction step if $p= m-1$ :
we consider the family of 
sections 
$$u^{k+1}\otimes s^{(0)}_j\in H^0(S, (k+1)L^{(m)}+ (m+k+1)A_{|S}),$$ 
In this case we have to exhibit a metric
on the bundle 
$$L_{m}+ L^{(m-1)}+kL^{(m)}+  (m+k+1)A.$$
This is easier than before, since we can simply take 
$$\varphi^{m(k+1)}:= 
 m\varphi_{\Xi}+ \varphi_A+ 
\varphi^{(km+m-1)}\leqno (33)$$
With this choice, the curvature conditions are satisfied ; as for the $L^2$ ones, we remark that we have
$$\eqalign{
 & \int_S{{|u^{k+1}\otimes s^{(0)}_j|^2}\over {(\sum_q \vert u^k\otimes s^{(m-1)}_q\vert^2)}}
\exp \big(-m\varphi_{\Xi}\big)dV< \cr
 < & C\int_S{{|u\otimes s^{(0)}_j|^2}\over {(\sum_q \vert s^{(m-1)}_q\vert^2)}}
\exp \big(-m\varphi_{\Xi}\big)dV< \infty \cr
}$$
where the last relation holds because of (29).
The proof of the extension theorem is therefore finished.\hfill \qed
\vskip 10pt

\bigskip
\noindent Next, we will highlight 
a very useful geometric context for which the 
numerous hypothesis in the theorem 2.1 are satisfied.

\noindent To this end, let us consider the next objects :
{\itemindent 5mm
\smallskip
\item {$(1)$} $\Delta_0$ is an effective $\bQ$-divisor and has critical exponent greater than 1
(i.e. it is klt) ;
\smallskip
\item {(2)} $\Delta_1, \Delta _2$ are effective, respectively ample $\bQ$--bundles ;
\smallskip
\item {(3)} The support of $S+ \Delta_0+ \Delta_1$ has normal crossing
and $S$ is not contained in the support of $\Delta_0+ \Delta_1$ ;
\smallskip
\item {(4)} Let $\Delta:= \Delta_0+ \Delta_1+\Delta_2$ ; we assume that there exist a closed positive current 
$$T\in c_1(K_X+ S+ \Delta)$$
such that the restriction of $T$ to any intersection of components of the support of the divisor $S+ \Delta_0+ \Delta_1$ is well defined (compare with the Hacon-McKernan
assumption on the stable base loci in [?]).

}

\medskip 
\noindent Our claim is that the assumptions in the theorem 2.1 are satisfied, 
provided that $(1)-(4)$ are. Indeed, let us first define the line bundles $L_j$
such that 
$$m\Delta_0= L_1+\cdots + L_{m-1}$$
and the corresponding metrics on them.

 By hypothesis, the $\bQ$--divisor $\Delta_0$ is effective, its critical exponent is
greater than one and the divisors in his support have normal crossing. Therefore, we have the decomposition
$$\Delta_0= \sum_{j= 1}^N{{p_j}\over {q_j}}Z_j$$
where the $(p_j, q_j)$ above are positive integers, such that coefficients above are {\sl strictly smaller   
than 1} ; we will assume that the sequence $\big({{p_j}\over {q_j}} \big)$ is increasing. 
As a consequence of the above decomposition, we see that for any positive and divisible enough 
$m$ we can find the line bundles 
$L_1,..., L_{m- 1}$ such that 
$$m\Delta_0\equiv L_1+...+ L_{m-1}$$
and such that 
$$L_j\equiv \sum_{p\in \Lambda_j} Z_p$$
with $\Lambda_j\subset \{1,...,N\}$. To see this, we first remark that 
for each index $j$ we have $\displaystyle m_j:= m{{p_j}\over {q_j}}\leq m-1$, and then 
we define 
$$\Lambda_j:= \{1,..., N\}$$ 
for $j= 1,...,m_1$. The next package of $(L_j)$ is defined 
by the set 
$$\Lambda_j:= \{2,...N\}$$
for $j:= m_1+1,..., m_2$, and so on. In this way we have our bundles 
$L_1,...L_{m-1}$~; we endow the corresponding $L_j$
with the natural singular metric induced by the hypersurfaces $(Z_p)$
whose indexes belong to $\Lambda_j$.  
We define 
$$L_m:= m(\Delta_1+ \Delta_3).$$
In view of the assumptions $(1)-(4)$ above,
all that still needs to be checked is the fact that the multiplier sheaves $(\cT_j)$ are trivial ; we are going to explain this along the following lines.

Let $y\in S$ be an arbitrary point, and let $j\in \{1,...,m-1\}$ be an arbitrary index.
We have to show that 
$$\int_{(S, y)}\exp \big((1-\delta)\varphi_{L_j}-\varepsilon \varphi_T\big)d\lambda< 
\infty$$
for all couple of parameters $0<\delta\ll \varepsilon\ll 1$.
In order to simplify the notations, we can assume that $L_j= Z_1+...+ Z_r$, and 
that $y\in S\cap Z_1\cap...\cap Z_b$ and $y\not\in Z_k$ for some $b\leq r, k\geq b+1$.
For each $p=1,..., b$ we define the complete intersection
$$S_p:= S\cap Z_1\cap...\cap Z_p$$
and the Skoda theorem [?] imply that 
$$\int_{(S_b, y)}\exp (-\varepsilon\varphi_T)d\lambda< \infty$$
for any $\varepsilon\ll 1$ (we remark that here we use the hypothesis $(4)$ 
concerning the restriction of $T$ to the sets $S_p$). 

By the local version of the theorem 2.2 above, we can extent the constant function equal to 1 on $S_b$ to a holomorphic function $f_{b-1}\in \cO(S_{b-1}, y)$ such that
$$\int_{(S_{b-1}, y)}|f_{b-1}|^2\exp \big(-(1-\delta)\varphi_{Z_{b}}- \varepsilon \varphi_T\big)d\lambda< \infty ;$$
we repeat this procedure $b$ times, until we get a function $f_0\in \cO(S, y)$
such that
$$\int_{(S, y)}|f_0|^2\exp \big(-(1-\delta)\varphi_{L_j}- \varepsilon \varphi_T\big)d\lambda< \infty.$$
Since the function $f_0$ is constant equal to 1 in a open set centered at $y$ in $S_b$,
we are done. \hfill\qed
\medskip
\noindent As a by-product of the preceding arguments, we have the following 
integrability criteria.

\claim Lemma|Let $\Theta$ be a closed (1,1)--current on a K\"ahler manifold
$(X, \omega)$, such that 
$$\Theta\geq -C\omega$$
for some positive constant $C$. Let $\displaystyle (Y_j)_{j\in J}$ be a finite set of hypersurfaces, which are supposed to be non-singular and to have normal crossings. Moreover, we assume that 
the restriction of $\Theta$ to the intersection $\cap _{i\in I}Y_i$ is well defined,
for any $I\subset J$. Then there exist a positive $\varepsilon_0= \varepsilon_0(\{\Theta\}, C)$ depending only on the cohomology class of the current and on the 
lower bound $C$ such that for any $\delta\in ]0, 1]$ and $\varepsilon \leq \varepsilon_0$ we have

$$\int_{(X, x)}\exp \big(-(1-\delta)\sum_{j\in J}\varphi_{Y_j}-\varepsilon\varphi_\Theta\big)d\Lambda< \infty,$$
at each point $x\in X$.
\endclaim
Indeed, the fact that the $\varepsilon_0$ only depends on the quantities in the above statement is a consequence of the fact that the Lelong numbers of closed positive currents on K\"ahler manifolds are bounded by the cohomology class of the current.
Also, it is very likely that the preceding statement hold true under the 
weaker assumption 
$$\nu_{\cap _{i\in I}Y_i}(\Theta)= 0$$
that is to say, we expect the previous lemma to be true if the generic Lelong number of $\Theta$ along all the intersections above is zero.
\hfill\qed

\bigskip

We remark next that the setting $(1)-(4)$ is precisely the context in which we have
found ourselves in the section 1.3. We recall now the result by Ein-Popa (see [14], as well as [7], [13], [26], [27], [28], [29], [31], [32], [33], [34]).

\claim 2.3 Theorem {\rm ([14])}|Let $X$ be a non-singular projective manifold. Let $S\subset X$ be a non-singular hypersurface,
and let $\Delta= \Delta_0+\Delta_1+\Delta_2$ be a sum of $\bQ$--line bundles on $X$
which verify the properties 1)-4) above.
Then any section $u$ of the bundle 
 $m(K_{X}+ S+ \Delta)_{|S}$ whose zero set contains the divisor $m\Delta_1$
 extends to $X$, for any $m$ divisible enough. \hfill \qed
 \endclaim
 
\proof. We will use the notations in the proof of the theorem 2.1 ; it was established
during this proof that the family of sections $\big(u^k\otimes s^{(0)}_j\big)$
 admits an extension to $X$, which was denoted by
 $$U^{(km)}_j\in H^0\bigl(X, kL^{(m)}+ mA\bigr)$$
 Here we remark that we have the factor $mA$ instead of $(k+m)A$
 in the 2.1 above, simply because $u$ is a section of the bundle 
 $m(K_{X}+ S+ \Delta)_{|S}$ (and not of $A+ m(K_{X}+ S+ \Delta)_{|S}$
 as in 2.1).
We take $k\gg 0$, so that
$$ {{m-1}\over {k}}A < \Delta_3,
$$
in the sense that the difference is
ample (we recall that $\Delta_3$ is ample). Let $h_A$ be a smooth, 
positively curved metric on $\Delta_3-{{m-1}\over{k}}A$. 
\vskip7pt
\noindent We apply now the extension theorem 2.2 with the following data
$$\displaystyle F:= {{m-1}\over {m}}L^{(m)}+ \Delta_1+ \Delta_2+ \Delta_3$$
and the metric $h$ on $F$ constructed as follows:
\smallskip
\item {(1)} on the factor  $\displaystyle {{m-1}\over {m}}L^{(m)}+ {{m-1}\over{k}}A$
we consider the appropriate power of the metric $\varphi^{(km)}$
given by the sections $(U^{(k)}_j)$; on $\Delta_1$ and $\Delta_2$ we consider the 
singular metric induced by the corresponding divisors, and finally on 
$ \Delta_3- {{m-1}\over{k}}A$ we take the metric $h_A$;
 
\smallskip
\item{(2)} we take an arbitrary, smooth metric $h_S$
on the bundle associated to the hypersurface $S$.
\smallskip
\noindent With our choice of the bundle $F$, 
the section $u$ we want to extend become a section of the adjoint 
bundle $K_X+ S+F_{\vert S}$ and the positivity requirements (1) and (2) in the extension
theorem are satisfied since we have
$$ {{\sqrt {-1}}\over {2\pi}}\Theta_h(F) \geq {{\sqrt {-1}}\over {2\pi}}\Theta_{h_A}(\Delta_3- {{m-1}\over
{mk}}B)> 0$$

Concerning the integrability of $u$,
remark that we have
$$\eqalign{
\int_S{{\vert u\vert^2\exp(- \varphi_{\Delta_1}- \varphi_{\Delta_2}-
\varphi_{A})}\over 
{\bigl(\sum_j\vert U^{(km)}\vert^2\bigr)^{{m-1}\over {mk}}}}= & 
\int_S{{\vert u\vert^2\exp(- \varphi_{\Delta_1}- \varphi_{\Delta_2}-
\varphi_{A})}\over 
{\bigl(\sum_j\vert u^{\otimes km}\otimes
s_j^{(0)}\vert^2\bigr)^{{m-1}\over {mk}}}}\leq \cr
\leq &C \int_S\vert u\vert^{2/m}\exp(-\varphi_{\Delta_2}-
\varphi_{A})\exp (- \varphi_{\Delta_1})dV<\cr
< & \infty\cr
}
$$
and this last integral converge by the hypothesis concerning the 
zero set of $u$ and the fact that the 
multiplier ideal sheaf of the restriction of $\Delta_1$ to $S$
is trivial.\hfill \qed
\medskip

\noindent 
\vskip 15pt

\vfill \eject 
\section{References}

\bigskip

{\eightpoint

\bibitem [1]&Berndtsson, B.:& On the Ohsawa-Takegoshi extension theorem;& Ann.\ Inst.\ Fourier (1996)&

\bibitem [2]&Berndtsson, B.:& Integral formulas and the Ohsawa-Takegoshi extension theorem;& Science in Chi\-na, Ser A Mathematics,  2005, Vol 48&

\bibitem [3]&Berndtsson, B.:& Curvature of Vector bundles associated to holomorphic fibrations;& to appear in Ann.\ of Maths.\ (2007)&

\bibitem [4]&Berndtsson, B., P\u aun, M.:& Bergman kernels and the pseudo-effectivity of the relative canonical bundles;& arXiv:math/0703344&

\bibitem [5]&P\u aun, M.:& A generalization of the Kawamata subadjunction theorem;& in  preparation&

\bibitem [6]&Boucksom, S.& Divisorial Zariski decompositions on compact complex manifolds;&  Ann. Sci. Ecole Norm. Sup. (4)  37  (2004),  no. 1, 45--76&

\bibitem [7]&Claudon, B.:& Invariance for multiples of the twisted canonical bundle~;& Ann. Inst. Fourier (Grenoble) 57 (2007), no. 1, 289--300&  

\bibitem [8]&Corti, A.:& Flips for 3-folds and 4-folds~;& Oxford University Press, Oxford, 2007&

\bibitem [9]&Demailly, J.-P.:& Regularization of closed positive currents and Intersection Theory; &J. Alg. Geom. 1 (1992) 361-409&

\bibitem [10]&Demailly, J.-P., Koll\'ar, J.:& Semicontinuity of complex singularity exponents and KŠhler-Einstein metrics on Fano orbifolds;& Ann. Ec. Norm. Sup 34 (2001), 525-556&

\bibitem [11]&Demailly, J.-P., Lazarsfeld, R.:& A subadditivity property of multiplier ideals;& Michigan Math. J. (special volume in honour of W. Fulton), 48 (2000), 137-156& 

\bibitem [12]&Demailly, J.-P.:&  On the Ohsawa-Takegoshi-Manivel  
extension theorem;& Proceedings of the Conference in honour of the 85th birthday of Pierre Lelong, 
Paris, September 1997, Progress in Mathematics, Birkauser, 1999&

\bibitem [13]&Demailly, J.-P.:& K\"ahler manifolds and transcendental techniques in algebraic geometry;&  Plenary talk and Proceedings of the Internat. Congress of Math., Madrid (2006), 34p, volume I&

\bibitem [14]&Ein, L., Popa, M.:& private communication;& june 2007\ &

\bibitem [15]&Favre, C., Jonsson, M.:&Valuations and plurisubharmonic functions&
J. Amer. math. Soc. 18 (205) 655-684&

\bibitem [16]&Hacon, C., D.,  McKernan, J.:& Boundedness of pluricanonical maps of varieties of general type;&
 Invent.\ Math.\ Volume {\bf 166}, Number 1 / October, 2006, 1-25&
 
\bibitem [17]&H\"oring, A.:& Positivity of direct image sheaves - a geometric point of view~;& in preparation, available on the author's home page&

\bibitem [18]&Kawamata, Y.:& Pluricanonical systems on minimal algebraic varieties~;& Invent. Math.  79  (1985),  no. 3&

\bibitem [19]&Kawamata, Y.:& On Fujita's freeness conjecture for
3-folds and 4-folds&  Math.\ Ann.\ {\bf 308}, 1997&

\bibitem [20]&Kawamata, Y.:& Subadjunction of log canonical divisors, 2;& Amer.\ J.\ Math.\  {\bf 120} (1998) 893--899&

\bibitem [21]&Lazarsfeld, R.:& Positivity in Algebraic Geometry;& Springer, Ergebnisse der Mathematik und ihrer Grenzgebiete&

 \bibitem [22]&Manivel, L.:& Un th\'eor\`eme de prolongement L2 de sections holomorphes d'un 
 fibr\'e hermitien;& Math.\ Zeitschrift {\bf 212} (1993), 107-122&

 \bibitem [23]&McNeal J., Varolin D.:&  Analytic inversion of adjunction: $L\sp
2$ extension theorems with gain;& Ann. Inst. Fourier (Grenoble) 57
(2007), no. 3, 703--718& 
  
 \bibitem [24]&Ohsawa, T., Takegoshi, K.\ :& On the extension of $L^2$
holomorphic functions;& Math.\ Z.,
{\bf 195} (1987), 197--204&
 
 \bibitem [25]&Ohsawa, T.\ :& On the extension of $L\sp 2$ holomorphic functions. VI. A limiting case;&
Contemp.\ Math. \ ,  {\bf 332} (2003)
Amer. Math. Soc., Providence, RI&

\bibitem [26]&P\u aun, M.:& Siu's Invariance of Plurigenera: a One-Tower Proof ;&preprint 2005, to appear in J.\ Diff.\ Geom.\ &

\bibitem [27]&Takayama, S:& On the Invariance and Lower Semi--Continuity
of Plurigenera of Algebraic Varieties;& J. Algebraic Geom.  {\bf 16 } (2007), no. 1, 1--18&

\bibitem [28]&Takayama, S:& Pluricanonical systems on algebraic varieties of general type;&Invent.\ Math.\
Volume {\bf 165}, Number 3 / September, 2005, 551-587&

\bibitem [29]&Tsuji, H.:& Extension of log pluricanonical forms from subvarieties;& math.CV/0511342 &
  
\bibitem [30]&Siu, Y.-T.:& Analyticity of sets associated to Lelong numbers and the extension of closed positive currents& Invent. Math.  27  (1974), 53--156&
 
\bibitem [31]&Siu, Y.-T.:& Invariance of Plurigenera;& Inv.\ Math.,
{\bf 134} (1998), 661-673&

\bibitem [32]&Siu, Y.-T.:& Extension of twisted pluricanonical sections with plurisubharmonic weight and invariance of semipositively twisted plurigenera for manifolds not necessarily of general type;& Complex geometry (G\"ottingen, {\bf 2000}),  223--277, Springer, Berlin, 2002&

\bibitem [33]&Siu, Y.-T.:& Multiplier ideal sheaves in complex and algebraic geometry;& Sci.\ China Ser.  {\bf A 48}, 2005&
 
\bibitem [34]&Varolin, D.:&  A Takayama-type extension theorem;&  math.CV/0607323, to appear in Comp.\ Math&
     
\bibitem [35]&Viehweg, E.:& Quasi-Projective
Moduli for
Polarized Manifolds;& Springer-Verlag, Berlin, Heidelberg, New York, 1995
as: Ergebnisse der Mathematik und ihrer Grenzgebiete, 3. Folge, Band 30&

}

\bigskip
\noindent
(version of April 23, 2008, printed on \today)
\bigskip\bigskip
{\parindent=0cm
Bo Berndtsson,  
bob@math.chalmers.se\\
Mihai P\u aun,
paun@iecn.u-nancy.fr
}

\end